\documentclass[11pt]{combine}
\usepackage{amsmath,amsfonts,amssymb,amsthm,epsfig, graphicx, hyperref}
\usepackage[dvipsnames,table,dvipsnames*, svgnames*, hyperref]{xcolor}
\usepackage{natbib,amsmath,amssymb,amsthm,graphicx,setspace,paralist,booktabs,rotating,subcaption,float,color}
\usepackage{multirow}
\usepackage{fancyhdr}
\usepackage{bibunits}
\usepackage[small]{titlesec}

\newtheorem{theorem}{Theorem}
\newtheorem{corollary}{Corollary}
\newtheorem{proposition}{Proposition}

\newtheorem{lemma}{Lemma}
{
	\theoremstyle{definition}
	\newtheorem{definition}{Definition}
	\newtheorem{example}{Example}

}
\newcommand{\beq}{\begin{equation}}
\newcommand{\eeq}{\end{equation}}
\newcommand{\beas}{\begin{align*}}
\newcommand{\eeas}{\end{align*}}
\newcommand{\bea}{\begin{align}}
\newcommand{\eea}{\end{align}}
\newcommand{\bei}{\begin{itemize}}
	\newcommand{\eei}{\end{itemize}}
\newcommand{\ben}{\begin{enumerate}}
	\newcommand{\een}{\end{enumerate}}
\newcommand{\bet}{\begin{theorem}}
	\newcommand{\eet}{\end{theorem}}
\newcommand{\bel}{\begin{lemma}}
	\newcommand{\eel}{\end{lemma}}
\newcommand{\bep}{\begin{proposition}}
	\newcommand{\eep}{\end{proposition}}
\newcommand{\bed}{\begin{definition}}
	\newcommand{\eed}{\end{definition}}
\newcommand{\bec}{\begin{corollary}}
	\newcommand{\eec}{\end{corollary}}
\newcommand{\bex}{\begin{example}}
	\newcommand{\eex}{\end{example}}

\newcommand{\Sig}{\bold{\Sigma}}

\newcommand{\bu}{\bold{u}}
\newcommand{\bw}{\bold{w}}
\newcommand{\bv}{\bold{v}}

\newcommand{\bfeta}{\boldsymbol{\eta}}
\newcommand{\balpha}{\boldsymbol{\alpha}}
\newcommand{\bbeta}{\boldsymbol{\beta}}

\newcommand{\R}{\mathbb{R}}
\newcommand{\E}{\mathbb{E}}

\newcommand{\calR}{\mathcal{R}}
\newcommand{\supp}{\text{supp}}

\newcommand{\argmax}{\mathop{\rm arg\max}}

\def\xx{\bold{x}}

\numberwithin{equation}{section}
\numberwithin{theorem}{section}
\numberwithin{lemma}{section}
\numberwithin{proposition}{section}
\numberwithin{remark}{section}

\begin{document}


\title{\scshape Optimal Estimation of Bacterial Growth Rates Based on Permuted Monotone Matrix}
\author{Rong Ma$^1$, T. Tony Cai$^2$ and Hongzhe Li$^1$ \\
Department of Biostatistics, Epidemiology and Informatics$^1$\\
Department of Statistics$^2$\\
University of Pennsylvania\\
Philadelphia, PA 19104}
\date{}
\maketitle
\thispagestyle{empty}

	\begin{abstract}
Motivated by the problem of estimating  the bacterial growth rates for genome assemblies from shotgun   metagenomic data, we consider the permuted monotone matrix model $Y=\Theta\Pi+Z$, where $Y\in \R^{n\times p}$ is observed, $\Theta\in \R^{n\times p}$ is an unknown approximately rank-one signal matrix with monotone rows, $\Pi \in \R^{p\times p}$ is an unknown permutation matrix, and $Z\in \R^{n\times p}$ is the noise matrix. This paper studies the estimation of the extreme values associated to the signal matrix $\Theta$, including its first and last columns, as well as their difference. Treating these estimation problems as compound decision problems, minimax rate-optimal estimators are constructed using the spectral column sorting method. Numerical experiments through simulated and  synthetic microbiome metagenomic data are presented, showing the superiority of the proposed methods over the alternatives. The methods are  illustrated by comparing  the growth rates of  gut bacteria between inflammatory bowel disease patients and  normal controls.  
	\bigskip
	
\end{abstract}

\section{Introduction}\label{intro.sec}

The statistical problem considered in this paper is motivated by  the problem of estimating the bacterial growth dynamics using shotgun metagenomics data. Several methods have been developed  to quantify the bacterial growth dynamics based on shotgun metagenomics data by extrapolating particular patterns in the sequencing read coverages  resulted from the bidirectional microbial DNA replications \citep{myhrvold2015distributed,abel2015sequence,korem2015growth, brown2016measurement}.  
For bacterial species with known complete genome sequences, \cite{korem2015growth} proposed to use the peak-to-trough ratio (PTR) of read coverages to quantify the bacterial growth rates after aligning the sequencing reads to the bacterial genomes.  Besides quantifying the growth rates for the bacteria with complete genome sequences, it is also of great importance  to estimate the growth rates of incomplete genome assemblies, where the coverages of contigs are observed in multiple samples. However, the order the contigs is only known up to an unknown permutation. 

Recently, \cite{gao2018quantifying} developed a computational algorithm (DEMIC) that accurately estimates the growth dynamics of a given assembled species by taking advantage of highly fragmented contigs assembled from multiple samples.  DEMIC is based on the following permuted monotone matrix model:
\beq \label{model}
Y=\Theta\Pi+Z
\eeq
where  the observed data $Y\in \R^{n\times p}$  is the matrix of  the preprocessed contig  coverage  for a given bacterial species. Specifically, the entry $Y_{ij}$ represents the log-transformed averaged read counts of the $j$-th contig of the bacterial species for the $i$-th sample after the preprocessing steps, including  genome assemblies, GC adjustment of read counts and outlier  filtering. In practice, the data set is usually high-dimensional in the sense that the number of contigs $p$  far exceeds the sample size $n$, so throughout we assume $p\gg n$. The signal matrix $\Theta\in \R^{n\times p}$  represents the true log-transformed coverage matrix of $n$ samples and $p$ contigs, where each row is  monotone due to the bi-directional DNA replication mechanism \citep{brown2016measurement,gao2018quantifying}, $Z\in \R^{n\times p}$ is the noise matrix, and $\Pi \in \R^{p\times p}$ is a permutation matrix, corresponding to some permutation $\pi$ from the symmetric group $\mathcal{S}_p$.  
\cite{ma2020optimal}  developed  methods for optimally recovering  the underlying permutation $\pi$ from $Y$. In particular, considering the loss function being either the 0-1 loss or the normalized Kendall's $\tau$ distance, a minimax optimal permutation estimator is proposed and theoretically analyzed under various parameter spaces.

\begin{figure}[h!]
	\centering
	\includegraphics[angle=0,width=4.5cm]{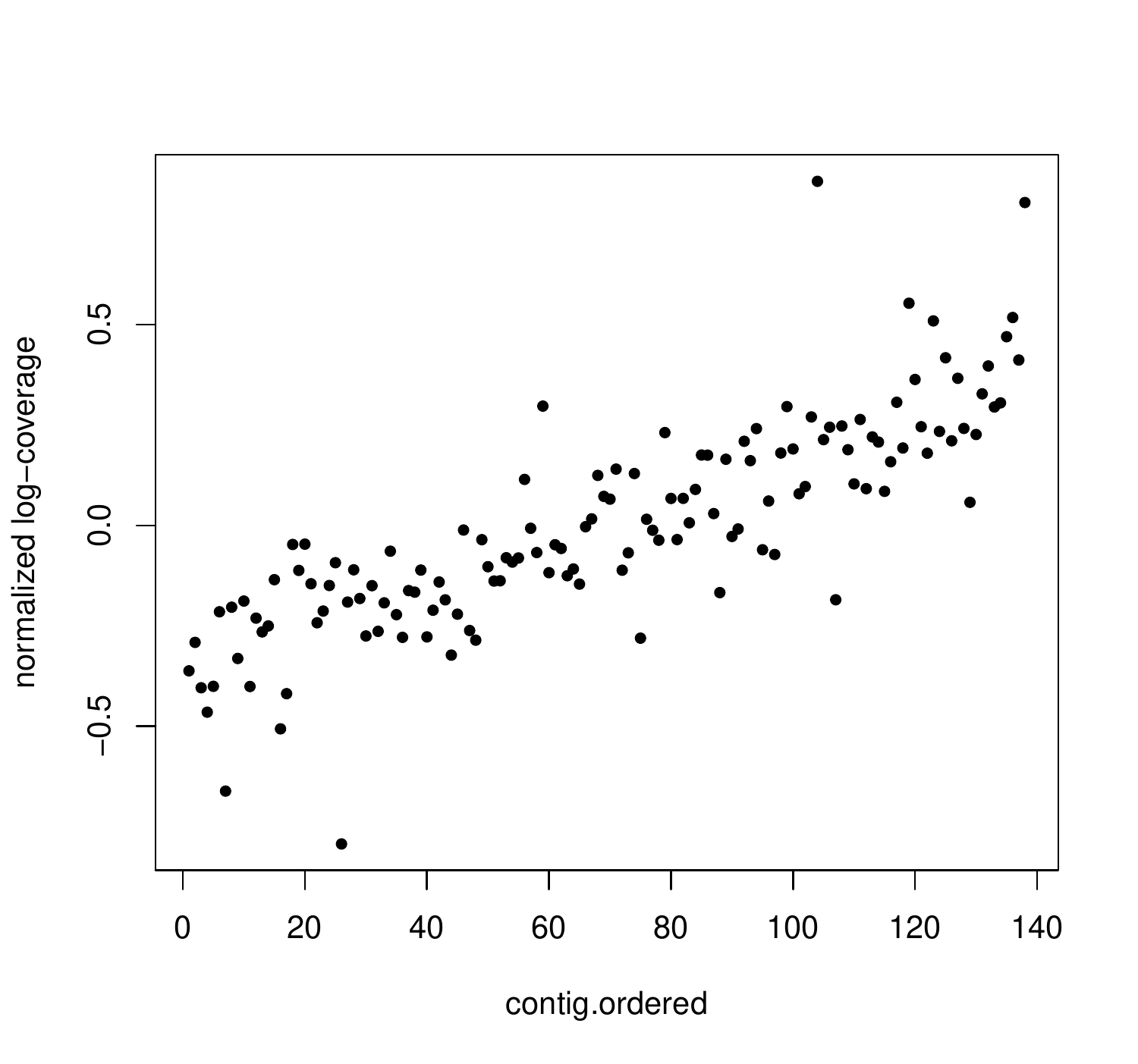}
	\includegraphics[angle=0,width=4.5cm]{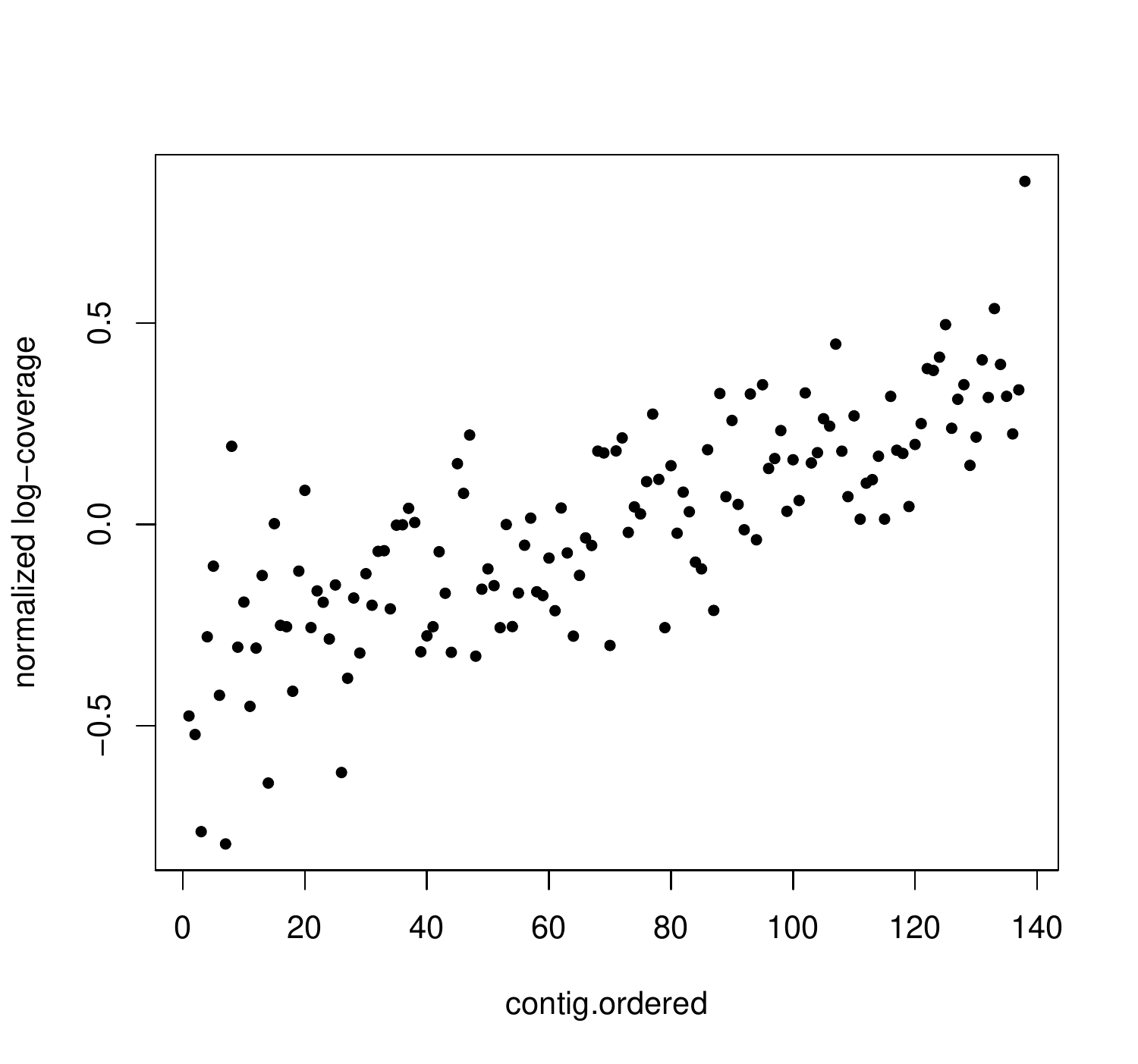}
	\includegraphics[angle=0,width=4.5cm]{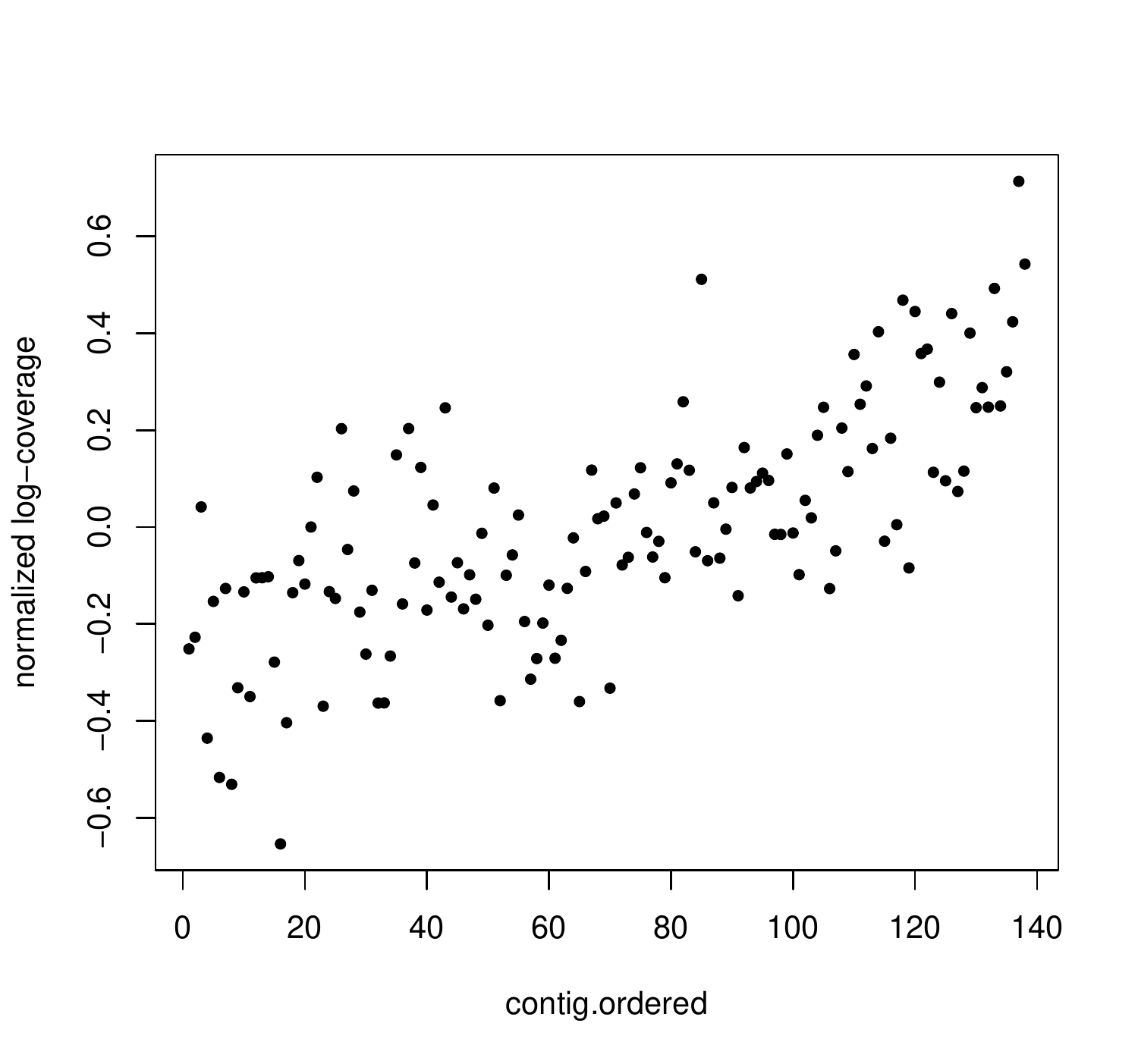}
	\caption{The log-coverages of ordered contigs of an assembled bacteria species   from 3 individuals with inflammatory bowel disease, detailed in Section \ref{real.data}.} 
	\label{evi.fig}
\end{figure}

In addition to the monotonicity constraint imposed on the rows of $\Theta$, real metagenomic data sets also suggest approximate linear relationship between the contig positions and their log-coverages for each sample, which indicates {approximately rank-one} structure of $\Theta$, after certain normalization. As an example, Figure \ref{evi.fig} shows the normalized log-contig counts of an assembled bacterial genome  for three individuals along estimated contig orders, suggesting the aforementioned approximate linear or rank-one structure, see Section \ref{risk.sec} for details. 

Under the permuted monotone matrix model, one can relate the two extreme columns $\Theta_R$ and $\Theta_L$, i.e., the first and the last columns of $\Theta$, to the log-transformed true peak and trough coverages of a given bacterial species, and define their difference  ${ R}(\Theta)=\Theta_R- \Theta_L$ as the true log-PTRs that characterize the bacterial growth rates  over $n$ samples.  
The goal of this paper  is to  provide a rigorous statistical framework for optimal estimation of the extreme values in the approximately rank-one permuted monotone matrix model,  including $\Theta_{R}$ and $\Theta_{L}$ and the range vector ${ R}(\Theta)$.
Based on the idea of spectral column sorting and the theory of low-rank matrix estimation, we develop computationally efficient estimators for the extreme columns and the range vector. In particular, 
the minimax optimality of the proposed methods are theoretically established and empirically illustrated with numerical experiments, which also justify its applicability in analyzing real data sets such as the microbiome metagenomics data.

Throughout the paper, we define the permutation $\pi$ as a bijection from the set $\{1,2,...,p\}$ onto itself. For simplicity, we denote $\pi=(\pi(1),\pi(2),...,\pi(p))$. All permutations of the set $\{1,2,...,p\}$ form a symmetric group, equipped with the function composition  operation $\circ$, denoted as $\mathcal{S}_p$. For any $\pi\in \mathcal{S}_p$, we denote $\pi^{-1}\in \mathcal{S}_p$ as its group inverse, so that $\pi\circ\pi^{-1}=\pi^{-1}\circ\pi=id$. In particular, we may use $\pi$ and its corresponding permutation matrix $\Pi\in \R^{p\times p}$ interchangeably, depending on the context. For a vector ${a} = (a_1,...,a_n)^\top \in \mathbb{R}^{n}$, we define the $\ell_p$ norm $\| {a} \|_p = \big(\sum_{i=1}^n a_i^p\big)^{1/p}$.  For a matrix $\Theta\in \R^{p_1\times p_2}$, we denote $\Theta_{.i}\in \R^{p_1}$ as its $i$-th column and denote $\Theta_{i.}\in \R^{p_2}$ as its $i$-th row. We write $a\land b=\min\{a,b\}$ and $a\lor b=\max\{a,b\}$. Furthermore, for sequences $\{a_n\}$ and $\{b_n\}$, we write $a_n = o(b_n)$ if $\lim_{n} a_n/b_n =0$, and write $a_n = O(b_n)$, $a_n\lesssim b_n$ or $b_n \gtrsim a_n$ if there exists a constant $C$ such that $a_n \le Cb_n$ for all $n$. We write $a_n\asymp b_n$ if $a_n \lesssim b_n$ and $a_n\gtrsim b_n$. 

\section{Extreme Value Estimation via Spectral Sorting} \label{extreme column}

\subsection{Spectral Sorting and Extreme Column Localization} \label{sorting}

A crucial step for estimating the extreme columns is to sort the permuted columns in order to identify the extreme ones. In this section, we introduce a spectral approach for localizing the permuted columns.
Toward this end, for any $\Theta$ with monotone rows, we consider the row-centered matrix
\beq \label{theta'}
\Theta'=\Theta \bigg( { I}_p-\frac{1}{p}{ ee}^\top \bigg)\in \R^{n\times p},
\eeq
where ${ e}=(1,....,1)^\top\in\R^p$.
Intuitively, $\Theta'$ is invariant to the row averages of $\Theta$ and preserves the row-monotonicity structure as well as the distances between the columns of $\Theta$.
The singular value decomposition of $\Theta'$ can be written as
\beq \label{Theta.svd}
\Theta'=\sum_{i=1}^r \lambda_i \bu_i\bv_i^\top,\qquad \text{for some $r\le \min\{n,p \}$},
\eeq
where $\lambda_1\ge \lambda_2\ge ...\ge \lambda_r$ are the ordered singular values of $\Theta'$ and $\bu_i$ and $\bv_i$ are the  left and right singular vectors corresponding to $\lambda_i$, respectively.
To overcome the identifiability issue, we assume \\
\indent ({A}) $\lambda_1$ has multiplicity one and the first nonzero component of $\bv_1$ is negative. \\
The following proposition provides an important insight that the row-monotonicity of a matrix actually implies the monotonicity of the components of its leading right singular vector $\bv_1$. This property plays a fundamental role in analyzing the permuted monotone matrix model.

\bep \label{rsv.prop}
Let $\Theta$ be a row-monotone matrix, whose row-centered version $\Theta'$ defined in (\ref{theta'}) satisfies (A).  Then its first right singular vector $\bv_1=(v_{11},...,v_{1p})^\top$ is a centered monotone vector, i.e., $\sum_{i=1}^p v_{1i}=0$ and $v_{11}\le v_{12}\le...\le v_{1p}$. In addition, the sign vector $\textup{sgn}(\bu_1)$ indicates the direction of monotonicity of the rows of $\Theta'$ (or $\Theta$).
\eep

From the above proposition, the relative orders of the columns of $\Theta'$ (and $\Theta$) are qualitatively preserved by the leading right singular vector $\bv_1$, whereas the directions of monotonicity for different rows are coded by the leading left singular vector $\bu_1$. As a result, given a column-permuted and noisy matrix $Y$ in (\ref{model}), one could localize the extreme columns $\Theta_R$ and $\Theta_L$ in $\Theta\Pi$ by considering the row-normalized observation matrix $X=Y({ I}_p-\frac{1}{p}{ ee}^\top)$ and its first right singular vector, i.e.,
\beq \label{hat.bv}
\widehat{\bv}=(\widehat{v}_1,...,\widehat{v}_p)^\top=\argmax_{\bv\in \R^p:\|v\|_2=1} \bv^\top X^\top X\bv.
\eeq
In accordance with Proposition \ref{rsv.prop}, it was shown by \cite{ma2020optimal} that the order statistics $\{\widehat{v}_{(1)},...,\widehat{v}_{(p)}\}$ can be used to optimally recover the permutation $\pi$, or the original column orders, by tracing back the permutation map between the elements of $\widehat{\bv}$ and their order statistics. Clearly, for  extreme column localization, the extreme values statistics $\widehat{v}_{(1)}$ and $\widehat{v}_{(p)}$ are more relevant. In fact, it is shown in the subsequent section  that, minimax optimal estimators can be constructed using such spectral extreme values estimates.

\subsection{Compound Decision Problem and the Proposed Estimators} \label{compound}

The problem of estimating $\Theta_R$, $\Theta_L$ or ${ R}$ consists of $n$ individual sub-problems, namely, estimating each of its $n$ coordinates.
Following the concept proposed by \cite{robbins1951asymptotically,robbins1964empirical} and further elaborated in \cite{samuel1967compound,copas1969compound, zhang2003compound} and \cite{brown2009nonparametric}, among many others, we observe that the problem of finding their minimax optimal estimators is a compound statistical decision problem, as the $n$ individual  sub-problems are amalgamated into one larger problem through the combined risk shown in equation  (\ref{risk}). Moreover, although the observations over $n$ samples are independent, it has been argued that, in general, for a compound decision problem, usually the simple estimators, where only the $i$-th sample is used to estimate the $i$-th coordinate, are suboptimal; in contrast, a minimax optimal estimator should be compound in the sense that multiple samples are used for the estimation of each coordinate.

In light of our discussion in Section \ref{sorting} as to the fundamental role of $(\lambda_1,\bu_1,\bv_1)$, we introduce our proposed  estimators for the extreme columns as
\beq \label{edge.est}
\widehat{\Theta}_{R}^*
=\widehat{v}_{(p)}X\widehat{\bv}+\frac{1}{p}Y { e} \in \R^n,\qquad
\widehat{\Theta}_{L}^*
=\widehat{v}_{(1)}X\widehat{\bv}+\frac{1}{p}Y { e} \in \R^n,
\eeq
and our proposed range estimator as
\beq
\widehat{ R}^*=\widehat{\Theta}_{R}^*-\widehat{\Theta}_{L}^*=(\widehat{v}_{(p)}-\widehat{v}_{(1)})X\widehat{v},
\eeq
where we recall that $\widehat{\bv}$ is defined in (\ref{hat.bv}) and $\widehat{v}_{(i)}$ is the $i$-th smallest order statistic among $\{\widehat{v}_1,...,\widehat{v}_p \}$.  By construction, the proposed extreme column estimators (\ref{edge.est}) are compound estimators, and each of them consists of two parts: the first part estimates the extreme columns of the row-centered matrix $\Theta'$ whereas the second part compensates the row-specific mean effects. In particular, in accordance with the observations made in Section \ref{intro.sec}, to construct the first parts of $\widehat{\Theta}_{R}^*$ and $\widehat{\Theta}_{L}^*$, the approximately rank-one structure $\Theta'_{.\ell}\approx\lambda_1v_{1\ell}\bu$ for $\ell\in\{1,p\}$, is incorporated with $v_{1\ell}$ estimated by $\widehat{v}_{(\ell)}$ and $\lambda_1 \bu_1$ estimated by $X\widehat{\bv}$.

\cite{ma2020optimal} developed an optimal estimator for the permutation $\pi$ under the model (\ref{model}).
Specifically,  let $\mathfrak{r}: \R^p \to \mathcal{S}_p$ be the ranking operator, which is defined such that for any vector $\xx\in \R^p$, $\mathfrak{r}(\xx)$ is the vector of ranks for components of $\xx$ in increasing order -- whenever there are ties, increasing orders are assigned from left to right.  The best linear projection estimator of $\pi$ was defined as $\hat{\pi}=[\mathfrak{r}({\widehat{\bv}})]^{-1}$.
This permutation estimator can be used to construct a natural two-step estimator of the two extreme columns. In the first step, we recover/sort the columns of $Y$ to obtain  the sorted matrix $\check{Y}=\begin{bmatrix} Y_{.\hat{\pi}(1)} & Y_{.\hat{\pi}(2)} &... & Y_{.\hat{\pi}(p)} \end{bmatrix}$. Intuitively, the column-sorted matrix $\check{Y}$ is expected to be close to $\Theta$.  In the second step, we fit a simple linear regression between each row of $\check{Y}$ and the sorted projection scores $(\widehat{v}_{(1)},\widehat{v}_{(2)},...,\widehat{v}_{(p)})$, which characterize the column relative locations. Denote the fitted intercepts as $\balpha=(\alpha_1,...,\alpha_n)^\top$ and the slopes as $\bbeta=(\beta_1,...,\beta_n)^\top$. We define the two-step regression estimators  as
\beq 
\hat{\Theta}_{L}^{Reg}=\balpha + \bbeta \widehat{v}_{(1)},\quad \hat{\Theta}_{R}^{Reg}=\balpha + \bbeta \widehat{v}_{(p)},\quad \hat{ R}^{Reg}=\bbeta( \widehat{v}_{(p)}-\widehat{v}_{(1)}).
\eeq
It is easy to check (see Section 3 of Supplementary Material) that under the conditions of Proposition \ref{rsv.prop}, it holds that
\beq \label{equiv}
\hat{\Theta}_{L}^{Reg}=\widehat{\Theta}_L^*,\qquad \hat{\Theta}_{R}^{Reg}=\widehat{\Theta}_R^*,\qquad \hat{ R}^{Reg}=\widehat{ R}^*.
\eeq
Intuitively, the extreme columns of the sorted matrix $\check{Y}$ could be suboptimal as it does not make use of the rank-one structure. A better way is to project rows of $\check{Y}$ onto the eigenspace spanned by $\widehat{ v}$, which is equivalent to regressing rows of $\check{Y}$ to $\widehat{ v}$.
This interesting observation provides another way of understanding  our proposed estimators.

\section{Theoretical Properties}
\subsection{Risk Upper Bounds for the Extreme Column Estimators} \label{risk.sec}

In what follows, we study the theoretical properties of our proposed estimator $\widehat{\Theta}_{R}^*$, as the results for $\Theta_L$ would hold in parallel.  Towards this end, we consider the normalized $\ell_2$ distances $\|\hat{\Theta}_{R}-\Theta_{R}\|_2/\surd{n}$ and denote the corresponding estimation risk as
\beq \label{risk}
\calR_{R}(\hat{\Theta}_{R})=\frac{1}{\surd{n}}E\|\hat{\Theta}_{R}-\Theta_{R}\|_2.
\eeq
We first define the set of monotone matrices
\[
\mathcal{D} = \Bigg\{ \Theta = (\theta_{ij})\in \R^{n\times p}: \quad \begin{aligned} 
& \text{for each $1\le i\le n$, either $\theta_{i,j}\le \theta_{i,j+1}$ for all $j$,}\\ 
&	 \text{or  $\theta_{i,j}\ge \theta_{i,j+1}$ for all $j$} 
\end{aligned} 
\Bigg\}.
\]
Recall that the row-centered version of $\Theta$, namely $\Theta'$, has the singular value decomposition given by (\ref{Theta.svd}). Consequently, throughout, we consider the following parameter space for $(\Theta,\pi)$
\begin{align} \label{D_R}
\mathcal{D}_{R}(t,\beta)&=\bigg\{ (\Theta,\pi)\in \mathcal{D}\times  \mathcal{S}_p: \begin{aligned}
&\text{(A) holds}, 0\le v_{1p}\le \beta,\\
&\lambda_1 \in[t/8,8t], \textstyle \sum_{i=2}^r\lambda_i\le \sigma\surd{\log p}
\end{aligned} \bigg\},
\end{align}
with $t\ge 0$ and $p^{-1/2}\le \beta\le 1$. Here the constraint on $\beta$ is natural since $\bv_1$ is a unit vector and $\beta$ is no less than the order of its largest component. Intuitively, the hyper-parameters $(t,\beta)$ characterize the global signal strength as well as the relative position of the extreme column $\Theta_R$ shared by the signal matrices in $\mathcal{D}_{R}(t,\beta)$, while the condition $\sum_{i=2}^r\lambda_i\le \sigma\surd{\log p}$ imposes a strong approximately rank-one structure on the row-centered $\Theta$. As our proposed estimators do not intend to estimate the possible additional structures upon the leading rank-one structure, such approximate rank-one condition is in this sense intrinsic to the problem.

To simplify notation, we define the rate function $\psi=\psi(n,p)=\surd{({\log p}/{n})}.$
The following theorem provides a uniform risk upper bound of the proposed estimator $\widehat{\Theta}_{R}^*$ over $\mathcal{D}_{R}(t,\beta)$.

\bet[Uniform Upper Bound] \label{unif.general.thm}
Suppose the pair $(t,\beta_R)$ satisfies $p^{-1/2}\le \beta_R\le 1$, $t^2\gtrsim \sigma^2\big[ \frac{1}{\beta_{R}^2} \land\big\{\frac{1}{\psi^2}+\frac{1}{\psi}\surd{(\frac{p}{n\log p})}\big\}\big]n\log p$ and the noise matrix $Z$ has independent sub-Gaussian entries $Z_{ij}$ with parameter $\sigma^2$. Then
\begin{align} \label{unibd.1}
\sup_{\mathcal{D}_{R}(t,\beta_R)}\calR_{R}(\widehat{\Theta}_{R}^*)&\lesssim \frac{\beta_{R}t}{\surd{n}}\bigg[\frac{\sigma\surd{\{(t^2+\sigma^2p)n\}}}{t^2}\land 1 \bigg]+\sigma\psi.
\end{align}
\eet

The risk upper bound (\ref{unibd.1})
consists of two components. In the first component, the factor $[{\sigma\surd{\{(t^2+\sigma^2p)n\}}}/{t^2}\land 1]$ is the error from estimating the leading left singular vector $\bu_1$ by its sample counterpart, whereas the factor ${\beta_R}{t}/\surd{n}$ reflects the overall magnitude of the extreme column $\Theta_R$ of the matrices in $\mathcal{D}_{R}(t,\beta_R)$. As for the second component $\sigma\psi(n,p)$, it comes from using the order statistic $\widehat{v}_{(p)}$ to estimate the largest component of $\bv_1$.

Interestingly, about the first component, we observe two phase transitions when $t^2$ passes $\sigma^2\surd{(np)}$ and $\sigma^2 p$, respectively. Specifically, in (\ref{unibd.1}), we have 
\begin{align*} 
\frac{\beta_{R}t}{\surd{n}}\bigg[\frac{\sigma\surd{\{(t^2+\sigma^2p)n\}}}{t^2}\land 1 \bigg]\asymp \left\{ \begin{array}{ll}
{\beta_Rt}/{\surd{n}}& \textrm{if $t^2\lesssim \sigma^2\surd{(np)}$,}\\
\frac{\beta_R\sigma^2\surd{p}}{t} & \textrm{if $ \sigma^2\surd{(np)}\lesssim t^2\lesssim \sigma^2p,$}\\
{\beta_R\sigma} & \textrm{if $t^2\gtrsim \sigma^2p$.}
\end{array} \right.
\end{align*}
From the theory of low-rank matrix estimation \citep{cai2018rate}, the quantity $\sigma^2\surd{(np)}$ is the critical point, below which it is impossible to estimate  the singular vector $\bu_1$. Hereafter we refer the collection of parameter spaces $\{\mathcal{D}_{R}(t,\beta_R):t^2\le \sigma^2\surd{(np)}\}$, $\{\mathcal{D}_{R}(t,\beta_R): \sigma^2\surd{(np)}\lesssim t^2\lesssim \sigma^2p\}$  and $\{\mathcal{D}_{R}(t,\beta_R): t^2\gtrsim\sigma^2p\}$  as the ``weak", ``  ``intermediate" and ``strong" signal-to-noise ratio regime, respectively.

To see the implications of the condition 
\beq \label{cond.1}
t^2\gtrsim \sigma^2\bigg[ \frac{1}{\beta_{R}^2} \land\bigg\{ \frac{1}{\psi^2}+\frac{1}{\psi}\surd{\bigg(\frac{p}{n\log p}\bigg)}\bigg\}\bigg]n\log p 
\eeq
of Theorem \ref{unif.general.thm}
on the critical events $t^2\asymp \sigma^2\surd{(np)}$ and $t^2\asymp \sigma^2 p$, we note  that, as long as $\beta_R\gg (n/p)^{1/4},$ by ignoring the logarithmic factors, the right-hand side of the condition (\ref{cond.1}) is asymptotically smaller than both critical points $\sigma^2\surd{(np)}$ and $\sigma^2 p$, so that both phase transitions exist under the condition of Theorem \ref{unif.general.thm}.

\subsection{Optimality of the Extreme Column Estimators and Minimax Rates}

Now we establish the minimax rate of convergence and the optimality of the proposed extreme column estimator $\widehat{\Theta}_R^*$ over the parameter space $\mathcal{D}_{R}(t,\beta_R)$. Specifically, for some given $(t,\beta)$, we define the minimax risks over $\mathcal{D}_{R}(t,\beta_R)$ as $\inf_{\hat{\Theta}_{R}}\sup_{\mathcal{D}_{R}(t,\beta_R)}\calR_{R}(\hat{\Theta}_{R})$ where the infimum is over all the possible estimators obtained from the data.
The following theorem provides the minimax lower bound of the estimation risk under the Gaussian noise.

\bet[Minimax Lower Bound] \label{homo.lower.thm}
Suppose $Z$ in model (\ref{model}) has i.i.d. entries $Z_{ij}\sim N(0,\sigma^2)$. Then, for any $\mathcal{D}_{R}(t,\beta)$ such that $t^2\ge c_0 \big( \frac{1-\beta_R^2}{\beta_R^2}\sigma^2\log p + \frac{\beta_R^2}{1-\beta_R^2}\sigma^2 p \big)$ and $c_1p^{-1/2}\surd{\log p}\le \beta_R\le c_2$, for sufficiently large $(n,p)$ and some constants $c_0,c_1>0$ and $0<c_2<1$ , it holds that
\begin{align}
\inf_{\widehat{\Theta}_R}\sup_{\mathcal{D}_{R}(t,\beta_R)} \calR_{R}(\widehat{\Theta}_{R})\gtrsim \frac{\beta_{R}t}{\surd{n}}\bigg[\frac{\sigma\surd{\{(t^2+\sigma^2p)n\}}}{t^2}\land 1 \bigg]+\sigma\psi. \label{ThetaR.lb1}
\end{align}
\eet

The proof of  Theorem \ref{homo.lower.thm} is involved. The main difficulty lies in the non-linearity and multi-dimensionality of the maps from the original parameter $\Theta$ to its extreme columns of interest. As the lower bound contains several components, we essentially derived three distinct minimax lower bounds corresponding to different worst-case scenarios. In addition to adopting the existing techniques such as the sphere packing of the Grassmannian manifolds,  we developed  two novel lower bound techniques to facilitate the proof of the minimax lower bound. The details can be found in the Supplementary Material.  The different conditions in Theorems 2 are due to the specific constructions in the lower bound argument of the proof. However, from a broader perspective, the conditions in Theorem 2 agrees to the ones in Theorem 1 in the sense that, the first condition, a lower bound on $t^2$, ensures that the global signal strength is sufficiently large, while the second condition is a mild restriction, up to a logarithmic factor, on the range of $\beta_R$. 

Combining the upper and the lower bounds, we obtain the exact minimax rate for estimating $\Theta_R$. Specifically, under the conditions of Theorems \ref{unif.general.thm} and  \ref{homo.lower.thm}, i.e., for $Z_{ij}\sim_{i.i.d.} N(0,\sigma^2)$ and
\beq \label{11}
t^2 \gtrsim \sigma^2\bigg[ \frac{1}{\beta_{R}^2} \land\bigg\{ \frac{1}{\psi^2}+\frac{1}{\psi}\surd{\bigg(\frac{p}{n\log p}\bigg)}\bigg\}\bigg]n\log p +\bigg( \frac{1-\beta_R^2}{\beta_R^2}\sigma^2\log p+ \frac{\beta_R^2\sigma^2p}{1-\beta_R^2}\bigg),
\eeq
we have
\begin{align} \label{adapt.c1}
\inf_{\widehat{\Theta}_R}\sup_{\mathcal{D}_{R}(t,\beta_R)}\calR_{R}(\widehat{\Theta}_{R})\asymp \frac{\beta_{R}t}{\surd{n}}\bigg[\frac{\sigma\surd{\{(t^2+\sigma^2p)n\}}}{t^2}\land 1 \bigg]+\sigma\psi,
\end{align}
where the optimal rate is attained by our proposed estimator $\widehat{\Theta}_R^*$.
To make better sense of  condition (\ref{11}), we note that, as long as $\beta_R\gtrsim (n/p)^{1/4}$, then by ignoring the logarithmic factors, condition (\ref{11}) is equivalent to $t^2 \gtrsim \sigma^2 \surd{(np)}$, which means the minimax rate can essentially be established over the intermediate to strong signal-to-noise ratio (SNR) regime, where the minimax rate is
\beq
\inf_{\widehat{\Theta}_R}\sup_{\mathcal{D}_{R}(t,\beta_R)}\calR_{R}(\widehat{\Theta}_{R})\asymp \frac{\sigma n^{1/4}\surd{(t^2+\sigma^2p)}}{p^{1/4}t}+\sigma\psi.
\eeq
As a consequence of the phase transition phenomena pointed out earlier, some interesting insights about the interplay between the global signal strength $t^2$, the dimensionality of the problem, the hardness of estimating $\Theta_R$ and that of estimating the leading left singular vector $\bu_1$, can be obtained. Specifically, we observe that (i) within the intermediate SNR regime ($\sigma^2\surd{(np)}\lesssim t^2\lesssim \sigma^2 p$),  increasing the signal strength $t^2$ will reduce the difficulty of estimating $\bu_1$ and therefore the rate for estimating $\Theta_R$, and (ii) within the strong SNR regime ($t^2\gtrsim \sigma^2 p$), the difficulty of estimating $\Theta_R$ no longer depends on $t^2$, as in this case the improved estimation of $\bu_1$ is neutralized by the increased magnitude of $\Theta_R$. Especially, all the above rate analysis is subjected to a possible lower bound of $\psi(n,p)$. 

Moreover, since the above minimax optimal rates are simultaneously attained by the proposed estimator  $\widehat{\Theta}_R^*$ regardless of the specific value of the underlying indices $(t,\beta_R)$, then, under the sub-Gaussian noise, $\widehat{\Theta}_R^*$ is minimax rate-adaptive over the collection of parameter spaces $\mathcal{C}=\{  \mathcal{D}_{R}(t,\beta_R): p^{-1/2}c_1\surd{\log p}\le\beta_R\le c_2<1,\text{ (\ref{11}) holds}\}$.
In particular, whenever $\beta_R\gtrsim (n/p)^{1/4}$, by ignoring the logarithmic factors, our proposed estimator is rate-optimally adaptive over the collection of parameter spaces lying in the intermediate to strong SNR regime, namely, $\mathcal{C}_{Adap}=\{  \mathcal{D}_{R}(t,\beta_R):
t^2\gtrsim \sigma^2\surd{(np)}\}$.

\subsection{Optimality of the Range Estimator and Minimax Rates} \label{range}

As a direct consequence of our previous results on the extreme column estimation, the theoretical properties of the range estimator $\widehat{R}^*$ can be obtained in the same manner.
Again we consider the normalized $\ell_2$ distances $\|\hat{ R}-{ R}(\Theta)\|_2/\surd{n}$ and denote the corresponding estimation risk as $\calR_{W}(\hat{ R})=\frac{1}{\surd{n}}E\|\hat{ R}-{ R}(\Theta)\|_2.$
Define the parameter space
\begin{align} \label{D_W}
\mathcal{D}_{W}(t,\beta_R,\beta_L)&=\bigg\{ (\Theta,\pi)\in \mathcal{D}\times \mathcal{S}_p:  \begin{aligned}
& \text{(A) holds}, \lambda_1\in[t/8,8t], \textstyle \sum_{i=2}^r\lambda_i\le \sigma\surd{\log p}, \\
&   -\beta_L\le v_{11}\le 0\le v_{1p}\le\beta_R,
\end{aligned}
\bigg\},
\end{align}
where $t\ge 0$, $p^{-1/2}\le \beta_R,\beta_L\le 1$, and define the function $q'(x,y,n,p)=\sigma^2n\log p \big[ \frac{1}{x^2} \land\big\{ \frac{1}{\psi^2}+\frac{1}{\psi}\surd{(\frac{p}{n\log p})}\big\}\big]+\big( \frac{1-x^2}{x^2}\sigma^2\log p+ \frac{y^2\sigma^2p}{1-y^2}\big).$
The following theorem establishes the minimax rate of convergence for estimating ${ R}(\Theta)$ and the minimax optimality and adaptivity of our proposed estimator $\widehat{ R}^*$.

\bet[Minimax Rates] \label{r.homo.upper.thm}
Let $\beta_W=\beta_R+\beta_L$. Suppose $t^2\ge c_0q'(\beta_R\land \beta_L,\beta_R\lor \beta_L,n,p)$,
$c_1p^{-1/2}\surd{\log p}\le \{\beta_R,\beta_L\} \le c_2$ for sufficiently large $(n,p)$ and some constants $c_0,c_1>0$ and $0<c_2<1$, 
and $Z$ has independent sub-Gaussian entries $Z_{ij}$ with parameter $\sigma^2$. Then 
\begin{align}
\inf_{\hat{ R}}\sup_{\mathcal{D}_{W}(t,\beta_R,\beta_L)} \calR_{W}(\hat{ R}) &\asymp \frac{\beta_{W}t}{\surd{n}}\bigg[\frac{\sigma\surd{\{(t^2+\sigma^2p)n\}}}{t^2}\land 1 \bigg]+\sigma\psi.
\end{align}
In particular, the minimax rates are simultaneously attained by the estimator $\widehat{ R}^*$.
\eet

\section{A Special Case:  Permuted Linear Growth Model}

In the previous sections, theoretical results are obtained for the general approximately rank-one matrices characterized by (\ref{D_R}) (\ref{D_W}) as well as the conditions of Theorems \ref{unif.general.thm} to \ref{r.homo.upper.thm}. One advantage is the rich row-monotonicity structures contained in such parameter spaces, which adapts well to real applications such as our motivating example in microbiome studies  where the noisy data sets are generated from the shotgun metagenomic sequencing, See \cite{boulund2018computational, gao2018quantifying} and Figure \ref{evi.fig}. However, in many cases such as classical theories of  the bacterial growth dynamics, an important subclass of the general permuted monotone matrix model has usually been considered for its heuristic simplicity and explanatory power. We refer this sub-model as the \emph{permuted linear growth model}, where (\ref{model}) holds over the restricted set
\begin{align*}
\mathcal{D}_{0}=\bigg\{ (\Theta,\pi)\in\mathcal{D}\times \mathcal{S}_p: \quad \begin{aligned}
&\theta_{ij}=a_i\eta_j+b_i,
\text{ where $a_i,b_i\in \R$ for $1\le i\le n$, }\\ 
&\text{$\eta_j\le \eta_{j+1}$ for $1\le j\le p-1$ and $\textstyle\sum_{j=1}^p\eta_j=0$.}
\end{aligned} \bigg\}.
\end{align*}
In other words, each row of $\Theta$ has a linear growth pattern with possibly different intercepts and slopes.  Denote ${ a}=(a_i)_{1\le i\le n}$, $\eta=(\eta_j)_{1\le j\le p}$ and ${ b}=(b_i)_{1\le i\le n}$. In this case, the parameters of interest have the expressions of $\Theta_R={ a}\eta_p$, $\Theta_L={ a}\eta_1$ and ${ R}={ a}(\eta_p-\eta_1)$. 

In the context of bacterial growth dynamics, the above model is commonly referred as the Cooper-Helmstetter model \citep{cooper1968chromosome,bremer1977examination} that  associates the copy number of genes with their relative distances to the replication origin.  Specifically, $a_i$ is the ratio of genome replication time and doubling time for the $i$th sample, $\eta_j$ is the distance from the replication origin for the $j$th contig, and $b_i$ is related to the read counts at the replication origin and the sequencing depth. Consequently, the extreme columns ${ a}\eta_p$ and ${ a}\eta_1$ correspond to the true log-transformed peak and trough coverages that are used to quantify the bacterial growth dynamics across the samples (see also Section \ref{simu.real} and \ref{real.data} for more details).

In the following, we discuss the consequences for the estimation of $\Theta_R$ under this special linear growth model, and the results for estimating $\Theta_L$ and ${ R}$ follow similarly.
By definition, the singular value decomposition (\ref{Theta.svd}) for $\Theta\in \mathcal{D}_0$ has a reduced form. Specifically, the row-centered matrix $\Theta'$ is exactly rank-one, where the leading right singular vector $\bv_1$ has components 
\beq \label{def11}
v_{1j}= \frac{\eta_j}{\|\eta\|_2},\quad \text{for $j=1,...,p,$} 
\eeq
and the largest singular value admits the expression
\beq \label{def12}
\lambda_1 =  \|{ a}\|_2 \|\eta\|_2.
\eeq
Intuitively, the set $\{v_{1j}\}_{1\le i<j\le p}$ characterize the exact normalized column positions of $\Theta'$ (and $\Theta$), while $\lambda_1$ summarizes the slope magnitude of the rows and the overall separateness of the columns.
Consequently, the risk upper bound obtained in Theorem \ref{unif.general.thm} has a reduced form, which has simpler and more intuitive interpretations. Specifically, for any given  $\Theta\in \mathcal{D}_R(t,\beta_R)$, we consider the following pointwise risk upper bound 
\beq \label{point.bnd}
\calR_{R}(\widehat{\Theta}_{R}^*)\lesssim \frac{v_{1p}\lambda_1(\Theta)}{\surd{n}}\bigg[\frac{\sigma\surd{\{(\lambda_1^2(\Theta)+\sigma^2p)n\}}}{\lambda_1^2(\Theta)}\land 1 \bigg]+\sigma\psi,
\eeq
induced by (\ref{unibd.1}) of Theorem  \ref{unif.general.thm}. 
With the reparametrizations (\ref{def11}) and (\ref{def12}), we can rewrite (\ref{point.bnd}) as
\beq\label{point.bnd2}
\calR_{R}(\widehat{\Theta}_{R}^*)\lesssim \frac{\eta_p\|{ a}\|_2 }{\surd{n}}\bigg[\frac{\sigma\surd{\{( \|{ a}\|_2^2 \|\eta\|_2^2+\sigma^2p)n\}}}{ \|{ a}\|_2^2 \|\eta\|_2^2}\land 1 \bigg]+\sigma\psi.
\eeq
Some observations about this risk upper bound are in order.

(a) Over the low SNR regime where $ \|{ a}\|_2^2 \|\eta\|_2^2\lesssim \sigma^2\surd{(np)}$, (\ref{point.bnd2}) becomes
\beq \label{l.3}
\calR_{R}(\widehat{\Theta}_{R}^*)\lesssim \frac{\|{ a}\|_2\eta_{p}}{\surd{n}}+\sigma\psi,
\eeq
where the first term is proportional to the overall slope magnitude $\|{ a}\|_2$, but does not rely on the locations of the other columns, i.e., $\eta_j$ for $1\le j\le p-1$. In this case, since the signal changes across different columns are so vague, $\widehat{\Theta}_R^*$ fails to implement a good estimate for the slopes ${ a}$ and the estimation error can only decrease when the extreme column $\Theta_R={ a}\eta_p$ itself (and its norm $\|{ a}\|_2\eta_p$) is close to zero. 

(b) Over the intermediate SNR regime where $ \sigma^2\surd{(np)}\lesssim  \|{ a}\|_2^2 \|\eta\|_2^2\lesssim \sigma^2 p$, (\ref{point.bnd2}) becomes
\beq \label{l.1}
\calR_{R}(\widehat{\Theta}_{R}^*)\lesssim \frac{\sigma\eta_p}{\|\eta\|_2}{\bigg(1+\frac{\sigma^2p}{\|{ a}\|_2^2\|\eta\|_2^2}\bigg)^{1/2}}+ \sigma\psi.
\eeq
In this case, as the signal differences between every consecutive columns are steep enough so that the slopes ${ a}$ can be well estimated, increasing $\|\eta\|_2$ or $\|{ a}\|_2$ would expand the advantage and therefore leads to a better estimate.

(c) Over the strong SNR regime where $ \|{ a}\|_2^2 \|\eta\|_2^2\gtrsim \sigma^2p$, the upper bound (\ref{point.bnd2}) becomes
\beq \label{l.2}
\calR_{R}(\widehat{\Theta}_{R}^*)\lesssim \frac{\sigma \eta_p}{\|\eta\|_2}+\sigma\psi.
\eeq
In the case, the advantage of large $\|{ a}\|_2$ has been exploited to extremity so that increasing $\|{ a}\|_2$ will no longer improve the performance of $\widehat{\Theta}_R^*$. 

Comparing the rates from (\ref{l.3}) to (\ref{l.2}), an interesting discrepant role played by the overall slope magnitude $\|{ a}\|_2$ can be observed. In general, the theoretical performance of $\widehat{\Theta}_R^*$ is clearly driven by the global SNR $\| { a}\|_2^2\|\eta\|_2^2/\sigma^2$, which measures the magnitude of the signal changes and the degree of monotonicity relative to the noise level. 

Following the same argument as the proof of Theorem \ref{homo.lower.thm}, the minimax optimality of our proposed estimator $\widehat{\Theta}_R^*$ can be also established under the permuted linear growth model. Specifically, if we define the indexed parameter space $
\mathcal{D}_{0,R}(t,\beta)=\{ (\Theta,\pi)\in \mathcal{D}_0: 0\le \eta_p/\|\eta\|_2\le\beta,  \|{ a}\|_2 \|\eta\|_2 \in [t/8,8t] \}$,
then it can be shown that for any pair $(t,\beta_R)$ such that (\ref{11}) holds and $p^{-1/2}\surd{\log p}\lesssim \beta_R\le c<1$, 
\begin{align*} 
\inf_{\widehat{\Theta}_R}\sup_{\mathcal{D}_{0,R}(t,\beta_R)}\calR_{R}(\widehat{\Theta}_{R})\asymp \frac{\beta_{R}t}{\surd{n}}\bigg[\frac{\sigma\surd{\{(t^2+\sigma^2p)n\}}}{t^2}\land 1 \bigg]+\sigma\psi,
\end{align*}
where the optimal rate is simultaneously attained by the proposed estimator $\widehat{\Theta}_R^*$.

\section{Numerical Studies} \label{numerical}

\subsection{Simulation with Model-Generated Data}

To demonstrate our theoretical results and to compare with alternative methods, we generate data from model (\ref{model}) with various configurations of the signal matrix $\Theta$.  Specifically, the signal matrix $\Theta=(\theta_{ij})\in \R^{n\times p}$ is generated under the following two regimes: (i) $S_1(n,p,\alpha)$: for any $1\le i\le n$, $\theta_{ij}=a_i \eta_j+ b_i$ for $1\le j\le p$, where $a_i\sim \text{Unif}(0,\alpha)$, $b_i\sim \text{Unif}(0,6)$ and $(\eta_1,...,\eta_p)=(-1,0,0,...,0,1)$; and (ii)  $S_2(n,p,\alpha)$: for any $1\le i\le n$, $\theta_{ij}=\log(1+a_i j+\beta_i)$ for $1\le j\le p$ where $a_i\sim \text{Unif}(0,\alpha)$ and $b_i\sim \text{Unif}(0,6)$.
By construction, $S_1(\alpha,n,p)$ belongs to the linear growth model whereas $S_2(\alpha,n,p)$ does not. The elements of $Z$ are drawn from i.i.d. standard normal distributions, and, without loss of generality, we set $\Pi={ I}_p$.

For the extreme column $\Theta_R$, we compare the empirical performance of our proposed estimator $\widehat{\Theta}_R^*$ with (i) the direct sorting estimator $\widetilde{ \Theta}_{R}$ defined as $\widetilde{\Theta}_R = Y_{.\hat{\pi}(p)}$, where $\hat{\pi}$ is defined in Section \ref{compound}; and (ii) the order statistic estimator $\check{\Theta}_{R}=(Y_{i,(p)})_{1\le i\le n}$, as all the rows of $\Theta$ are monotonic increasing. For the range vector ${ R}(\Theta)$,  we compare our proposed estimator $\widehat{ R}^*$ with (i) the direct sorting estimator  $\widetilde{ R}_{DS}=Y_{.\hat{\pi}(p)}-Y_{.\hat{\pi}(1)}$, and (ii) the order statistic estimator $\widetilde{ R}_{OS}=(Y_{i,(p)}-Y_{i,(1)})_{1\le i\le n}$.
We use the empirical risk, or the averaged normalized $\ell_2$ distance, to compare these methods. 
For each setting, we evaluate the empirical performance of each method over a range of $n$, $p$ and $\alpha$. Each setting is repeated for 200 times. 

\begin{figure}[!ht]
	\centering
	\includegraphics[angle=0,width=13cm]{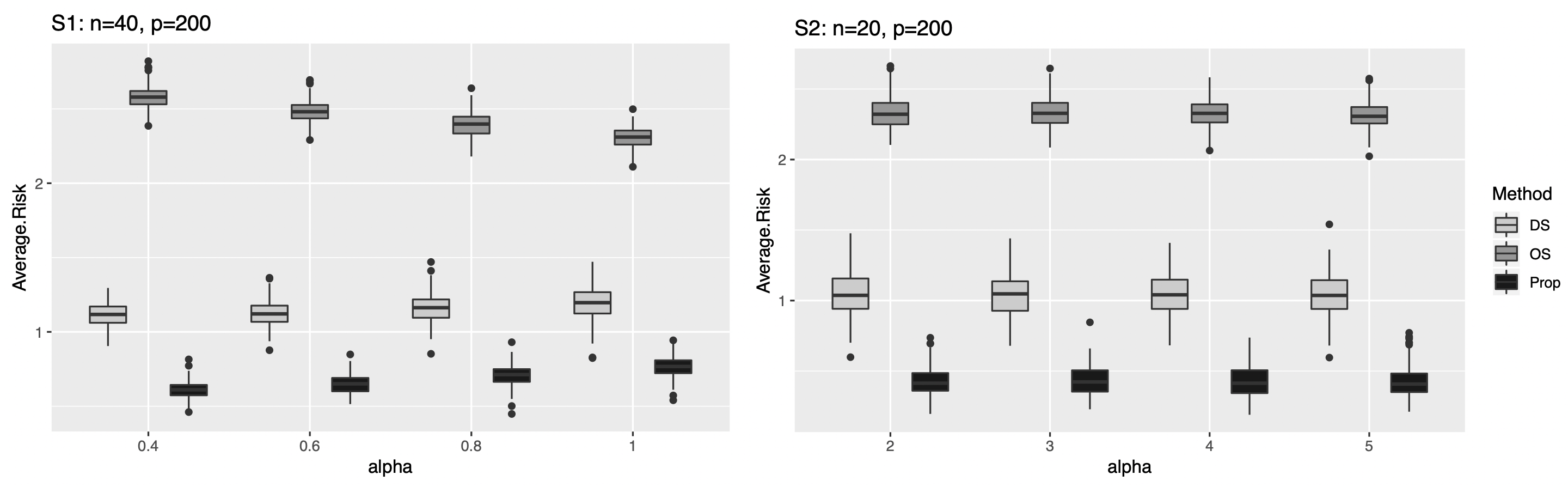}
	\includegraphics[angle=0,width=13cm]{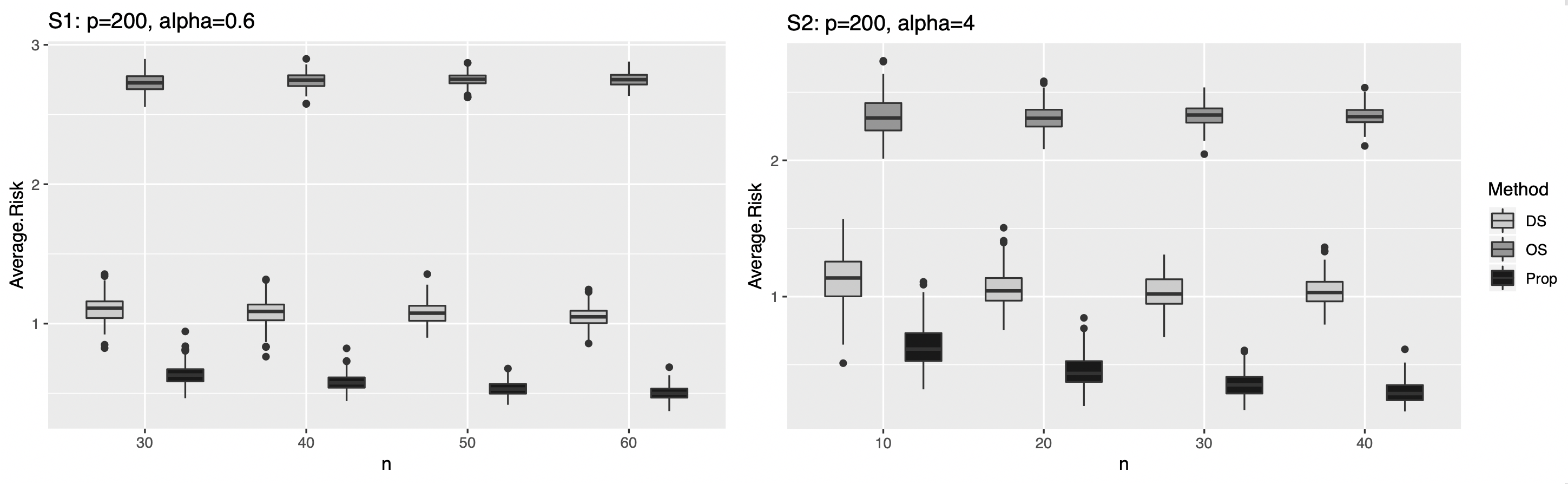}
	\includegraphics[angle=0,width=13cm]{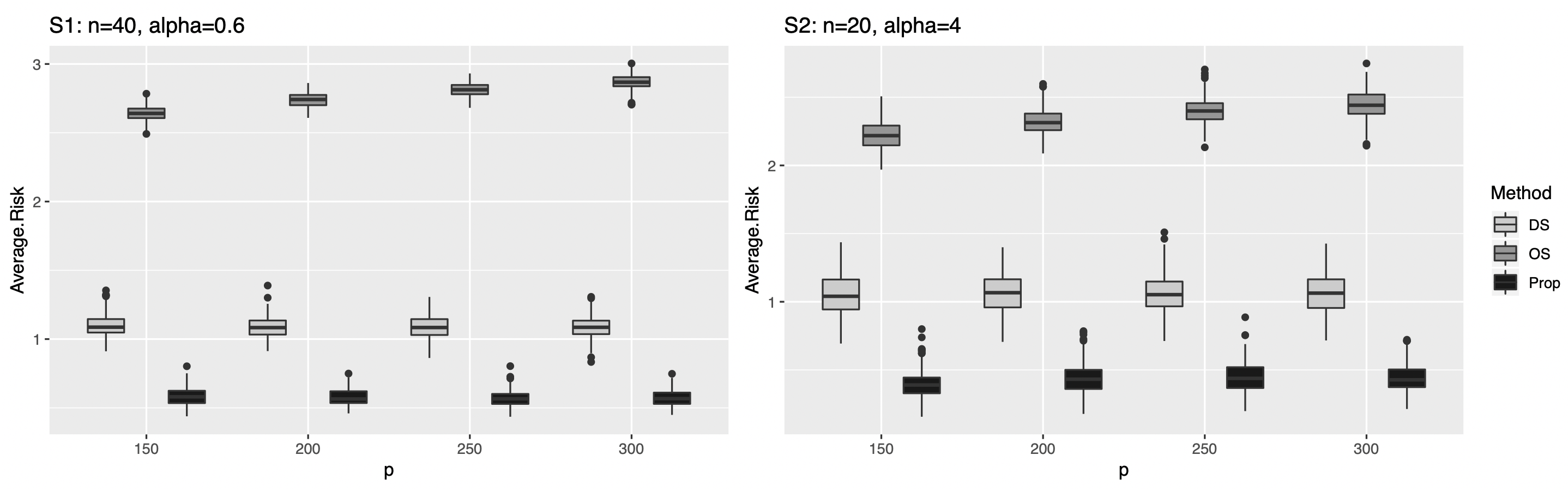}
	\caption{Boxplots of the empirical risks for estimating $\Theta_R$, with DS, OS and Prop representing $\widetilde{\Theta}_{R}$, $\check{\Theta}_{R}$ and $\widehat{\Theta}^*_R$.} 
	\label{simu.1}
\end{figure}

\begin{figure}[!ht]
	\centering
	\includegraphics[angle=0,width=13cm]{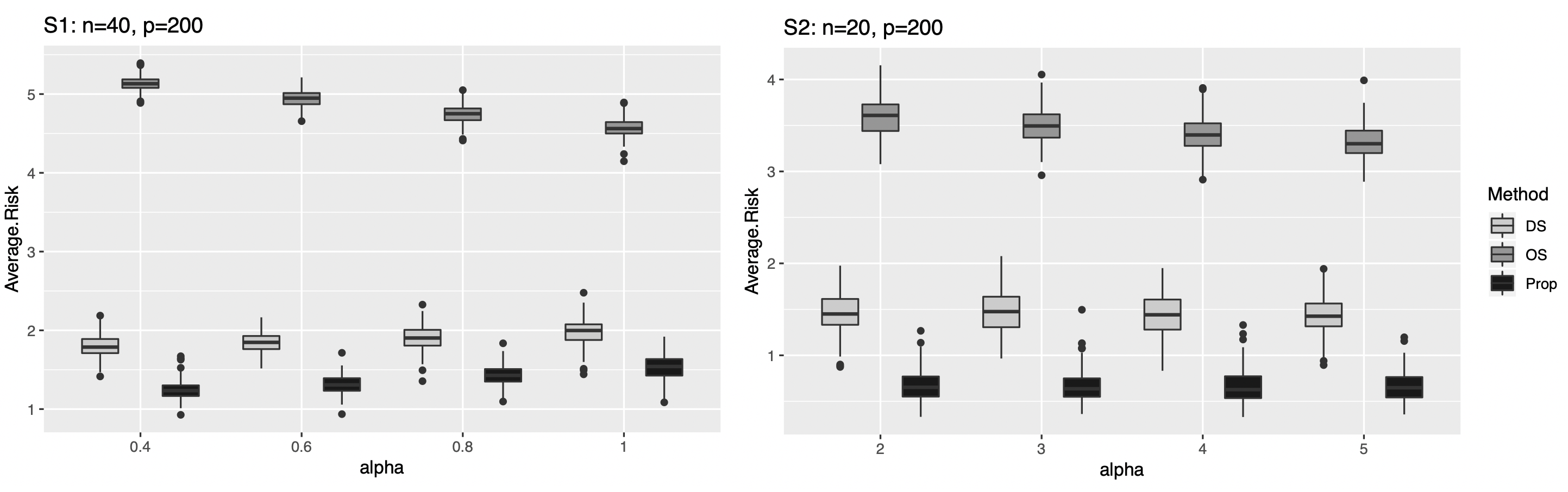}
	\includegraphics[angle=0,width=13cm]{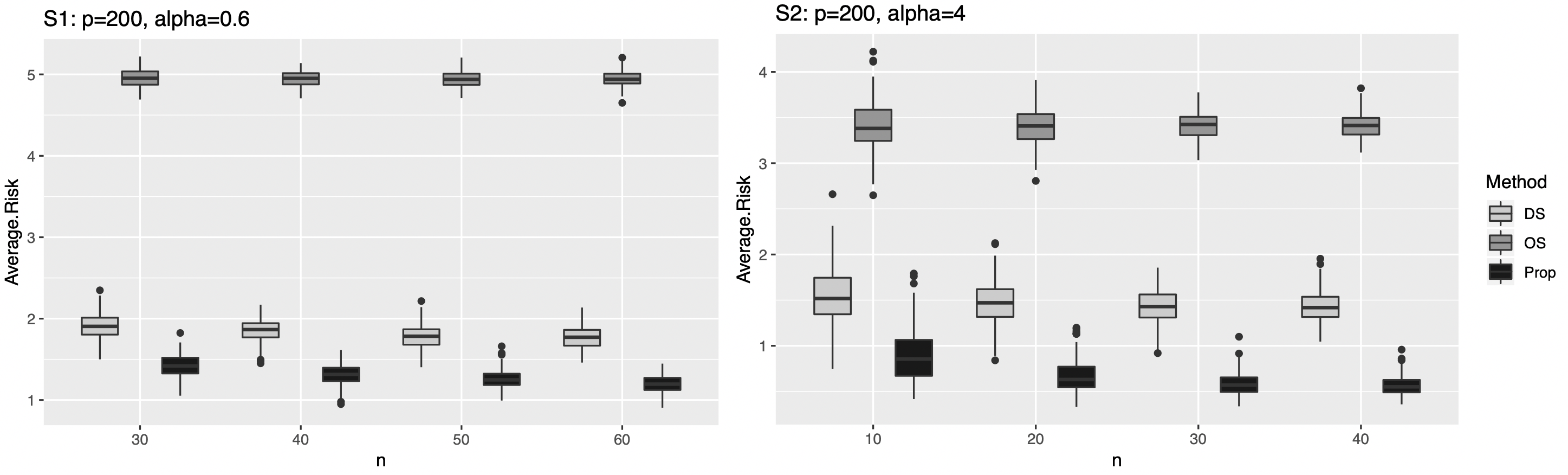}
	\includegraphics[angle=0,width=13cm]{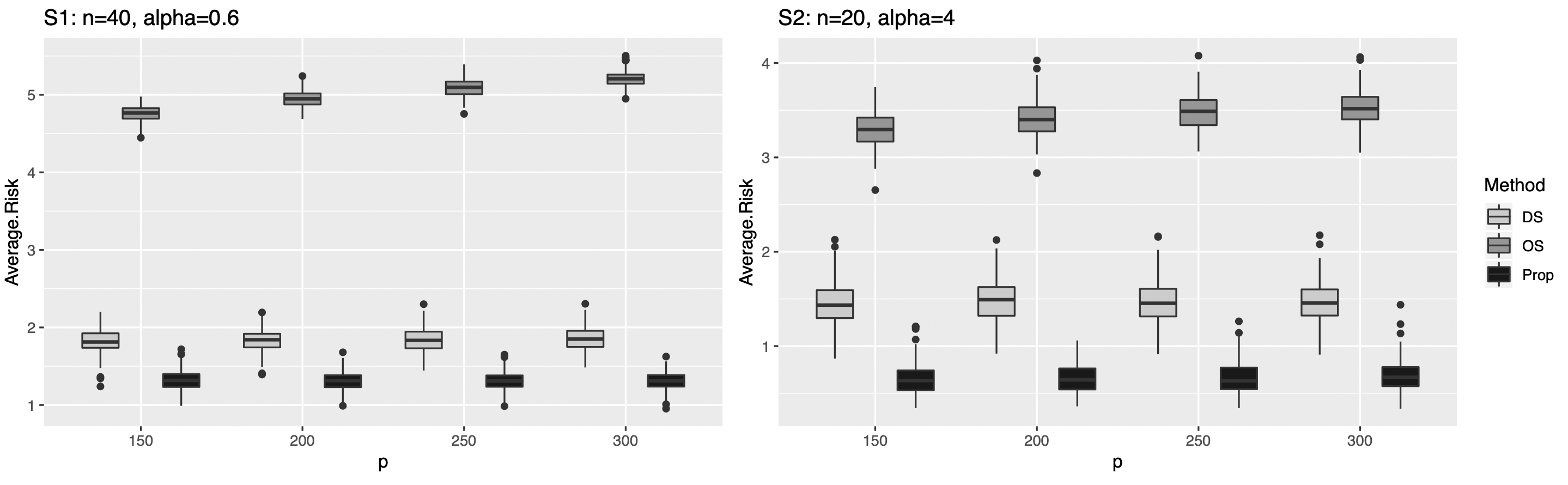}
	\caption{Boxplots of the empirical risks for estimating ${ R}$, with DS, OS and Prop representing $\widetilde{ R}_{DS}$, $\widetilde{ R}_{OS}$ and $\widehat{ R}^*$.} 
	\label{simu.2}
\end{figure}

The results are summarized as boxplots in Figure \ref{simu.1} and Figure \ref{simu.2}. The empirical results agree with our theory in the following perspectives: (i) our proposed estimators $\widehat{\Theta}_R^*$ and $\widehat{ R}^*$ perform the best among all the settings; (ii) in the middle two plots of Figure \ref{simu.1} and \ref{simu.2}, the risks of our proposed estimators decrease as $n$ grows, which agrees with our theorems. In addition, in the top left panel of Figure  \ref{simu.1} and \ref{simu.2} we observe that the risks of the order statistic estimator decrease as $\alpha$ increases. This is because under $S_1(\alpha,p,n)$, the parameter $\alpha$ characterizes the separateness of the two extreme columns from the other columns. The order statistic estimators would apparently favour the cases where the separation is more significant. In addition, both our proposed estimators and the direct sorting estimators outperform the order statistic estimators, showing the advantage of the compound estimators.

\subsection{Simulation with Synthetic Microbiome Metagenomic Data} \label{simu.real}

We then evaluate the empirical performance of our proposed method using a synthetic metagenomic sequencing data set   \citep{gao2018quantifying} by generating sequencing reads based on 45 closely related bacterial genomes in 50 independent samples. 
Particularly, 
\cite{gao2018quantifying} presented a synthetic shotgun metagenomic sequencing data set of a community of 45 phylogenetically related species from 15 genera of five different phyla with known RefSeq ID, taxonomy and replication origin \citep{synth}. To generate the metagenomic reads, reference genome sequences of three randomly selected  species in each genus were downloaded.  Read coverages were generated along the genome based on an exponential distribution with a specified peak-to-trough ratio and a function of accumulative distribution of read coverages along the genome was calculated. Sequencing reads were then generated using the above accumulative distribution functions and a random location for each read on the genome, until the total read number achieved a randomly assigned average coverage between 0$\cdot$5 and 10 folds for the species in a sample. Sequencing errors including substitution, insertion and deletion were simulated in a position- and nucleotide-specific pattern according to the metagenomic sequencing error profile of Illumina. 

For the final data set, the average nucleotide identities between species within each genus ranged from 66$\cdot$6\% to 91$\cdot$2\%.  The probability of one species existing in each of the 50 simulated samples was set as 0$\cdot$6, and a total of 1,336 average coverages and the corresponding PTRs were randomly and independently assigned.  After the same processing, filtering, and CG-adjustment steps as  in \cite{gao2018quantifying}, the final data set included genome assemblies of 41 species. For each species, we obtained the permuted matrix of log-contig coverage with the number of samples ranging from 29 to 46 and the number of contigs from 47 to 482. 

We provide estimates of the log-PTRs of the assembled species for all the samples, or the range vector ${ R}$, using our previous notations.   As a comparison, in addition to our proposed method $\widehat{ R}^*$, we consider the iRep estimator proposed by \cite{brown2016measurement}, where the contigs of a given species were ordered for each sample separately based on the observed read counts, before fitting a piece-wise linear regression function. We evaluate these methods by considering the $\ell_2$ distance between the vectors of the true log-PTRs and their estimates. To generalize our evaluation to diverse metagenomic data sets, we also evaluate the effect of sample size $n$ as well as contig numbers $p$ by randomly selecting subsets of samples or contigs from each data set. The selection was made with replacement. 

\begin{figure}[!ht]
	\centering
	\includegraphics[angle=0,width=14cm]{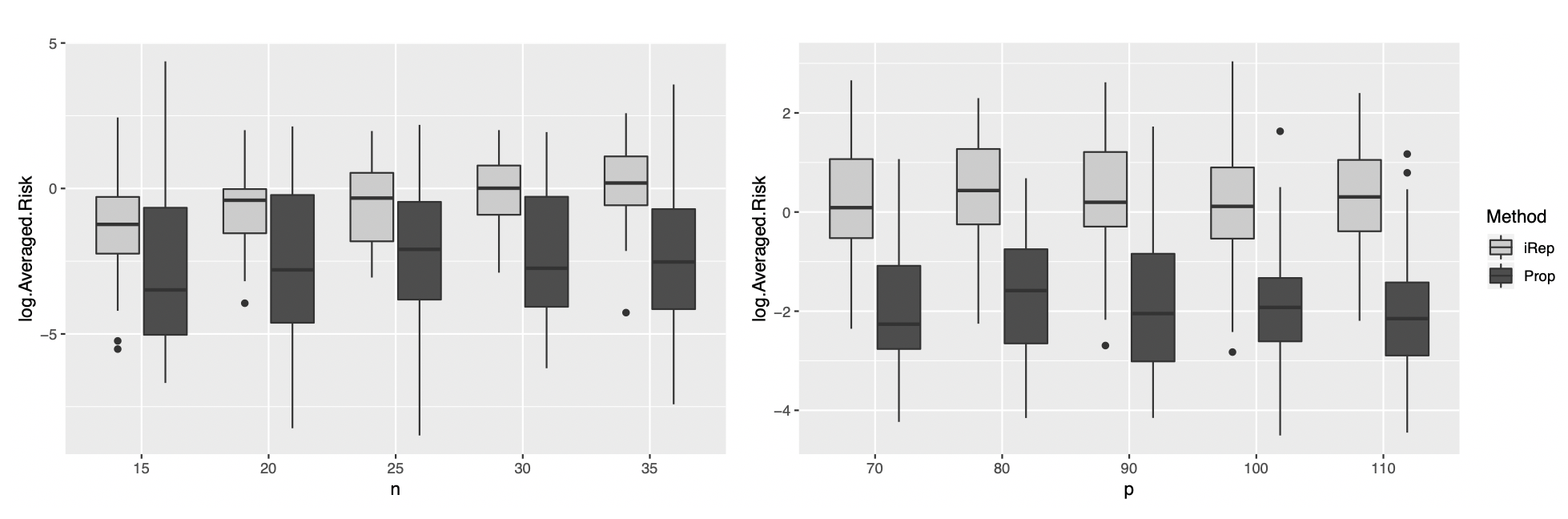}
	\caption{Boxplots of the $\ell_2$ distances between the estimated and the true log-PTRs. Prop: the proposed method; iRep: the iRep estimation method.}
	\label{synth.data}
\end{figure}

The results are summarized in Figure \ref{synth.data}.  As $n$ or $p$ varies, our proposed estimator performs consistently better than iRep. Moreover, the performance of our proposed method is not sensitive to the sample size, the number of contigs from the genome assemblies or the underlying true PTRs. These results partially explain why the DEMIC algorithm has superior performance compared to the existing ones \citep{gao2018quantifying}.

\subsection{Analysis of A Real Microbiome Metagenomic Data Set} \label{real.data}

We complete our numerical study by analyzing a real metagenomic data set from the NIH Integrative Human Microbiome Project. The project includes  the Inflammatory Bowel Disease (IBD) Multi'omics project to investigate the differences in  gut microbiome communities among adults and children with IBD \citep{lloyd2019multi} and normal non-IBD controls. Many studies  have reported strong associations between IBD, including  both Crohn's disease (CD) and ulcerative colitis (UC) and gut microbiota composition. In contrast, we focus on comparing the bacterial growth rates between UC, CD and normal non-IBD  individuals using the proposed methods.  

\begin{table}
	\centering
	\caption{\label{tab01} Analysis of bacterial growth rates among CD, UC and non-IBD samples. Bins that show significantly different growth rates and their taxonomic annotations are presented. $(n_1, n_2, n_3)$: numbers of CD, UC and non-IBD samples that carried the respective bin. } 
	\begin{tabular}{lccl}
		Bins &  Taxonomic Annotations & ($n_1,n_2,n_3$) & P-values  \\  
			\hline  
		bin$\cdot$054&    Roseburia (genus)& (54, 32, 54) & 0$\cdot$015 \\
		bin$\cdot$090 &   Faecalibacterium (genus)&(38, 41, 52) & 0$\cdot$005  \\
		bin$\cdot$091 &     Clostridiales (order)&(26, 40, 52) & 0$\cdot$016  \\
		bin$\cdot$099 & Subdoligranulum (genus)&(30, 32, 49) & $<0$$\cdot$$001$ \\
		bin$\cdot$465 & Dialister (genus)&(36, 41, 33)& 0$\cdot$043 \\
	\end{tabular}
\end{table}

The metagenomic data sets, including 300 samples of the CD, UC and non-IBD  subjects, were downloaded from the IBDMDB website, \texttt{https://www.ibdmdb.org}. Specifically, we randomly select 100 samples of UC, CD and normal non-IBD samples, respectively.  For each sample,  the sequencing data was obtained from the stool sample using Illumina shotgun sequencing.  We first apply  MEGAHIT \citep{li2015megahit} version 1$\cdot$1$\cdot$1  to perform  metagenomic co-assembly. The co-assembled contigs were then clustered into  metagenomic bins  or genome assemblies using MaxBin \citep{wu2015maxbin} version 2$\cdot$2$\cdot$4. Finally, Bowtie 2 \citep{langmead2012fast} version 2$\cdot$3$\cdot$2 was used to align reads back to the assembled contigs for each of the samples, and the output alignments were then sorted by samtools \citep{li2009sequence} version 0$\cdot$1$\cdot$19. 

After these preparations, the DEMIC algorithm, incorporated with our proposed methods, was applied to obtain the estimated PTRs (ePTRs) of a given species represented by a contig cluster (bin) for each sample. As a result, ePTRs of 25 bins were obtained for  subsets of the UC ($n_1$), CD ($n_2$) and non-IBD ($n_3$) samples with $n_1+n_2+n_3\ge 100$, as some contig clusters may not be carried or abundant enough among many samples.  For each bin, we compare the ePTRs among the  UC, CD and non-IBD samples using an F-test.  We applied the CAT/BAT algorithm \citep{von2019robust} that compares the metagenomic assembled bins to a taxonomy database  to obtain the taxonomic annotations of the 25 bins.  We observe that  only a few bins can be annotated at the species level, whereas many of the bins can  only be annotated to  genera or orders, suggesting that many of the assembled contig bins may correspond to new species. This agrees with a recent paper \citep{Almeida2020} that showed that more than 70\% of the assembled genomes lack cultured representatives.

Interestingly, based on the F-test, among the 25 contig clusters, 5 of them show significant difference in ePTRs among the UC, CD and non-IBD samples (Table \ref{tab01}). For reasons of space,  Table 1 only provides the taxonomic annotation of the bins in terms of their genus -- except for bin$\cdot$091 which can only be determined up to orders (see our Supplementary Material   for the complete annotations). We also performed pairwise comparisons using two-sample t-test  for the 5 differential bins (Table \ref{tab02}). We  found that the difference in the growth rates of bin$\cdot$054, Roseburia, bin$\cdot$090, Faecalibacterium, and bin$\cdot$099, Subdoligranulum, are more significant between IBD and non-IBD samples. In particular, boxplots in our Supplementary Material indicate higher growth rates of bin$\cdot$054, Roseburia,  and bin$\cdot$090, Faecalibacterium, and a lower growth rate of bin$\cdot$099, Subdoligranulum for IBD samples when compared  to the non-IBD samples. Moreover, the growth rates of bin$\cdot$091, Clostridiales, is significantly higher among UC samples, whereas the growth rate of bin$\cdot$465, Dialister, is significantly higher among the CD samples, comparing to the samples of the other two categories. These results show that the gut microbiome communities in CD and UC patients or IBD and non-IBD patients differ not only in relative abundance but also in growth rates of certain bacterial species, an important insight from our data analysis. 

\begin{table}
	\centering
	\caption{\label{tab02} p-values from pairwise $t$-tests of differential growth rates between different groups  for five genome assembly bins. }	
	\begin{tabular}{lcccc}
		Bins &  Taxonomic Annotations &UC vs. CD  &  UC vs. non-IBD & CD vs. non-IBD \\  
		\hline  
		bin$\cdot$054&    Roseburia (genus)&   0$\cdot$525 &  0$\cdot$004  & 0$\cdot$081\\
		bin$\cdot$090  &   Faecalibacterium (genus)&  0$\cdot$392  &  0$\cdot$016 & 0$\cdot$004\\
		bin$\cdot$091&     Clostridiales (order) & 0$\cdot$012 &   0$\cdot$054 & 0$\cdot$335 \\
		bin$\cdot$099 & Subdoligranulum (genus)& 0$\cdot$960 & $<0$$\cdot$$001$   & $<0$$\cdot$$001$\\
		bin$\cdot$465& Dialister (genus) & 0$\cdot$042& 0$\cdot$818  & 0$\cdot$026\\
	\end{tabular}
	
\end{table}


\section{Discussion} \label{dis}

The present paper focused on the permuted monotone matrix model with homoskedastic noise. If the noises are heteroskedastic, for example (i) the columns of the noise matrix are not independent, or (ii) the variances of the noise matrix entries are not identical, we argue that, as long as the marginal distributions of the noise matrix remain sub-Gaussian, the framework developed in this paper can still be applied. Specifically, in light of the recent work of \cite{zhang2018heteroskedastic}, where heteroskedastic principal component analysis and singular value decomposition are  studied, the key analytical tools paralleling to those applied in the current work, such as concentration and perturbation inequalities associated with the heteroskedastic random matrices, can be obtained by generalizing the results of \cite{zhang2018heteroskedastic}. Such extensions are involved and we leave them for future research.

The current theoretical framework was developed upon the approximately rank-one structure suggested by our specific metagenomic applications. However, extensions to other settings are possible by modifying the proposed methods.  In particular, the key observations made in Section \ref{sorting} apply to any monotone matrix satisfying condition (A). When the approximate rank-one assumption is violated, say, if the monotone signal matrix is rank-$r$ with $r>1$, one could construct  estimators based on the leading $r$ singular values and singular vectors by following the same idea in Section \ref{compound}, although its theoretical analysis might be technically challenging. 

\section*{Acknowledgements}

We would like to thank the Editor, Associate Editor and two anonymous referees for their constructive suggestions that improved the presentation of the paper. The research reported in this publication was supported by  grants from the National Institutes of Health and National Science Foundation. R. M. would like to thank Rui Duan, Yuan Gao and Shulei Wang for stimulating discussions at various stages of this project.

\section*{Supplementary materials}

Supplementary material includes the detailed proofs of the main theorems, and some supplementary tables and figures.

\bibliography{reference}
\bibliographystyle{chicago}

\newpage

\title{Supplement to ``Optimal Estimation of Bacterial Growth Rates Based on Permuted Monotone Matrix"}
\author{Rong Ma$^1$, T. Tony Cai$^2$ and Hongzhe Li$^1$ \\
	Department of Biostatistics, Epidemiology and Informatics$^1$\\
	Department of Statistics$^2$\\
	University of Pennsylvania\\
	Philadelphia, PA 19104}
\date{}
\maketitle
\thispagestyle{empty}

\begin{abstract}
In this Supplementary Material we prove the theorems and propositions in our main paper as well as the technical results. Some supplementary tables and figures are also included.
\end{abstract}

\setcounter{section}{0}

\setcounter{proposition}{1}

\section{Optimal Extreme Columns Estimation}

\subsection{Risk Upper Bound: Proof of Theorem \ref{unif.general.thm}}

It suffices to obtain the pointwise risk upper bound for a given $\Theta\in \mathcal{D}_R(t,\beta_R)$ and take the supremum over the parameter space.

First observe that the proposed estimator as well as its estimation risk is invariant to the unknown permutation matrix $\Pi$. Specifically, as both the calculation of $\widehat{\Theta}^*_{R}$ and the definition of $\Theta_{R}$ doesn't require knowledge about the permutation matrix $\Pi$, the quantity $\E \|\widehat{\Theta}^*_{R}-\Theta_{R}\|_2$ is invariant to the underlying permutation $\pi\in \mathcal{S}_p$. In the following, without loss of generality, we assume $\Pi={\bf I}_p$. 

Observe that 
\[
\widehat{\Theta}_R^*=\widehat{v}_{(p)}X\widehat{\bv}+\frac{1}{p}Y{\bf e}=\widehat{\lambda} \widehat{v}_{(p)} \widehat{\bu}+\frac{1}{p}Y{\bf e},\quad \Theta_R=\sum_{i=1}^r \lambda_i \bu_i v_{ip}+{\bf b,}
\]
where $\widehat{\lambda}$ and $\widehat{\bu}\in \R^n$ are first singular value and first left singular vector of $X$, $v_{ip}$ is the $p$-th component of $\bv_i$, and ${\bf b}=\frac{1}{p}\Theta {\bf e}$. We have
\begin{align}
\|\widehat{\Theta}_{R}^*-\Theta_{R}\|_2&\le  \|\widehat{\lambda}\widehat{\bu}\widehat{v}_{(p)}-\lambda_1 \bu_1 v_{1p}\|_2+\sum_{i=2}^r\|\lambda_i\bu_iv_{ip}\|_2+\bigg\|\frac{1}{p}Ye-{\bf b}\bigg\|_2.
\end{align}
As a result, as long as $\sum_{i=2}^r\|\lambda_i\bu_iv_{ip}\|_2=\sum_{i=2}^r\lambda_i|v_{ip}|\lesssim \|\widehat{\lambda}\widehat{u}\widehat{v}_{(p)}-\lambda_1 \bu_1 v_{1p}\|_2$,
it suffices to obtain upper bounds for $\E \|\widehat{\lambda}\widehat{\bu}\widehat{v}_{(p)}-\lambda_1 \bu_1 v_p\|_2$ and $\E \big\|\frac{1}{p}Y{\bf e}-{\bf b}\big\|_2$. In particular, we will show that
\beq \label{homo.upper.eq1}
\E  \|\widehat{\lambda}\widehat{\bu}\widehat{v}_{(p)}-\lambda_1 {\bu}_1 v_p\|_2 \lesssim \sigma\sqrt{\log p}+(\lambda_1 v_{1p}+\sigma \sqrt{n\log p}) \bigg(\frac{\sigma\sqrt{(\lambda_1^2+\sigma^2p)n}}{\lambda_1^2-\lambda_2^2}\land 1 \bigg),
\eeq
and
\beq \label{homo.upper.eq2}
\E \big\|\frac{1}{p}Y{\bf e}-{\bf b}\big\|_2 \lesssim \sigma\sqrt{\frac{n}{p}}.
\eeq
These two upper bounds along with the condition $\sum_{i=2}^r \lambda_i\le \sigma\sqrt{\log p}$ will lead to the risk bound
\beq
\calR_{R}(\widehat{\Theta}_{R}^*)\lesssim\bigg(\frac{\lambda_1v_{1p}}{\sqrt{n}}+\sigma\sqrt{\log p}\bigg)\bigg(\frac{\sigma\sqrt{(\lambda_1^2+\sigma^2p)n}}{\lambda_1^2-\lambda_2^2}\land 1 \bigg)+\sigma\psi,
\eeq
which implies
\beq \label{32}
\calR_{R}(\widehat{\Theta}_{R}^*)\lesssim \frac{v_{1p}\lambda_1}{\sqrt{n}}\bigg(\frac{\sigma\sqrt{(\lambda_1^2+\sigma^2p)n}}{\lambda_1^2}\land 1 \bigg)+\sigma\psi.
\eeq
if $\lambda_1\gg\lambda_2$ and either
\[
\sigma\sqrt{\log p}\lesssim \frac{\lambda_1v_{1p}}{\sqrt{n}}, \quad\text{ or }\quad \sigma\sqrt{\log p}\frac{\sigma\sqrt{(\lambda_1^2+\sigma^2p)n}}{\lambda_1^2-\lambda_2^2}\lesssim \sigma\psi.
\]
It can be checked that the above inequalities can be implied by
\beq \label{q.cond0}
\lambda_1^2\gtrsim \sigma^2n\log p \bigg[\frac{1}{v_{1p}^2} \land \bigg( \frac{1}{\psi^2}+\frac{1}{\psi}\sqrt{\frac{p}{n\log p}}\bigg)\bigg],
\eeq
which also implies $t^2\gtrsim \sigma^2\log p$, or $\lambda_1^2 \asymp \lambda_1^2-\lambda_2^2$. Thus we have (\ref{32}), which leads to the final uniform risk bound.
The rest of the proof is devoted to the upper bounds (\ref{homo.upper.eq1}) and (\ref{homo.upper.eq2}).

\paragraph{Proof of (\ref{homo.upper.eq1}).} Throughout, we write $v_{1p}$ as $v_p$ for convenience. By the elementary inequality that, for $x_i,y_i\ge 0, i=1,2,$
\begin{align}\label{basic.ineq}
|x_1x_2-y_1y_2|&\le y_1|x_2-y_2|+y_2|x_1-y_1|+|x_1-y_1||x_2-y_2|,
\end{align}
we have
\begin{align} \label{homo.upper.decomp}
\E\|\widehat{\lambda}\widehat{\bu}\widehat{v}_{(p)}-\lambda_1 {\bu}_1 v_p\|_2 &\le \E|\widehat{\lambda}\widehat{v}_{(p)}-\lambda_1 v_p|+\lambda_1 v_p\E\|\widehat{\bu}-{\bu}_1\|_2+\E[ |\widehat{\lambda}\widehat{v}_{(p)}-\lambda_1 v_p|\|\widehat{\bu}-{\bu}_1\|_2] \nonumber \\
&\lesssim \E|\widehat{\lambda}\widehat{v}_{(p)}-\lambda_1 v_p|+\lambda_1 v_p\E\|\widehat{\bu}-{\bu}_1\|_2,
\end{align}
where the last inequality follows from the trivial bound $\|\widehat{\bu}-{\bu}_1\|_2\le \sqrt{2}$.
It suffices to obtain upper bounds for $\E \|\widehat{\bu}-{\bu}_1\|_2$ and $\E  |\widehat{\lambda}\widehat{v}_{(p)}-\lambda_1 v_p|$.

Note that $\|\widehat{\bu}-{\bu}_1\|_2^2=2(1-\widehat{\bu}^\top {\bu}_1).$
By definition, up to a change of sign for $\widehat{\bu}$, we have $0\le \widehat{\bu}^\top {\bu}_1\le 1$. Then $1-\widehat{\bu}^\top {\bu}_1\le 1-(\widehat{\bu}^\top {\bu}_1)^2$ and
\beq \label{u.ineq}
\E\|\widehat{\bu}-{\bu}_1\|_2^2\le \E |1-(\widehat{\bu}^\top {\bu}_1)^2|.
\eeq
The upper bound can be obtained from the relation $\E \|\widehat{\bu}-{\bu}_1\|_2 \le \sqrt{\E \|\widehat{\bu}-{\bu}_1\|_2^2}$,
and the following lemma, which is a direct consequence of Proposition 1 and Theorem 3 of \cite{cai2018rate}.

\bel[\cite{cai2018rate}] \label{pca.lem2}
Let $X=\Theta+Z \in \R^{n\times p}$, where  $\Theta$ has SVD $\Theta=\sum_{i=1}^r \lambda_i u_i v_i^\top$ for $r\le \min\{n,p\}$ and $\lambda_1\ge \lambda_2\ge ...\ge \lambda_r$, $Z$ has independent sub-Gaussian entries  with parameter $\sigma^2$, and $\widehat{u}, \widehat{v}$ is the first left and right singular vectors of $X$, respectively. Then it follows that $\E |1-(\widehat{u}^\top u_1)^2| \le \frac{C\sigma^2(\lambda_1^2+\sigma^2p)n}{(\lambda_1^2-\lambda_2^2)^2}\land 1$ and $\E |1-(\widehat{v}^\top v_1)^2| \le \frac{C\sigma^2(\lambda_1^2+\sigma^2n)p}{(\lambda_1^2-\lambda_2^2)^2}\land 1$.
\eel

Applying the above lemma, 
\beq \label{uhat-u}
\E \|\widehat{\bu}-{\bu}_1\|_2\lesssim \bigg(\frac{\sigma\sqrt{(\lambda_1^2+\sigma^2p)n}}{\lambda_1^2-\lambda_2^2}\land 1 \bigg).
\eeq
Next, we show that
\beq \label{homo.max.bnd}
\E |\widehat{\lambda}\widehat{v}_{(p)}-\lambda_1 v_p| \lesssim \sigma\sqrt{\log p}+(\lambda_1 v_{1p}+\sigma \sqrt{n\log p}) \bigg(\frac{\sigma\sqrt{(\lambda_1^2+\sigma^2p)n}}{\lambda_1^2-\lambda_2^2}\land 1 \bigg).
\eeq
In particular, we need the following lemma concerning an elementary inequality.

\bel \label{max.ineq.lem}
Let $\{ a_i \}_{1\le i\le n}$ and $\{ b_i\}_{1\le i\le n}$ be two sets of real numbers. It then holds that $\big|\max_{1\le i\le n} a_i\big|-\big|\max_{1\le i\le n} b_i\big| \le \max_{1\le i\le n}|a_i-b_i|$.
\eel

Observe that $\widehat{\lambda}\widehat{v}_{(p)}= \max_{1\le j\le p}\widehat{\bu}^\top X_j$,
where $X_j$ is the $j$-th column of $X$. 
Using the inequality of Lemma \ref{max.ineq.lem}, we have
\begin{align} \label{scaler.bnd}
\E|\widehat{\lambda}\widehat{v}_{(p)}-\lambda_1 v_p| &= \E | \max_{1\le j\le p}(\widehat{\bu}^\top X_j-{\bf u}_1^\top X_j+{\bf u}_1^\top X_j-\lambda_1 v_p)| \nonumber\\
&\le \E \max_{1\le j\le p}|(\widehat{\bu}-{\bf u}_1)^\top X_j|+\E |\max_{1\le j\le p}{\bf u}_1^\top X_j-\lambda v_p| \nonumber \\
&\le \E \max_{1\le j\le p}|(\widehat{\bu}-{\bf u}_1)^\top \Theta'_j|+ \E \max_{1\le j\le p}|(\widehat{\bu}-{\bf u}_1)^\top E_j|+\E |\max_{1\le j\le p}{\bf u}_1^\top X_j-\lambda_1 v_p|
\end{align}
where we used the notation $X=\Theta'+E$ with $E_{ij}$ i.i.d. from $N(0,\sigma^2(p-1)/p)$.
The first term in the last inequality can be bounded by H\"older's inequality, the relation
\[
\max_{1\le j\le p}\|\Theta'_j\|_2=\max_{1\le j\le p} \bigg\| \sum_{i=1}^r\lambda_i  {\bu}_iv_{ij}  \bigg\|_2\le \lambda_1v_p+\max_{1\le j\le p} \sum_{i=2}^r\lambda_i\le \lambda_1v_p+\sigma\sqrt{\log p},
\] 
and (\ref{uhat-u}) as
\beq \label{homo.max.1}
\E \max_{1\le j\le p}|(\widehat{\bu}-{\bu}_1)^\top \Theta'_j|\le \E  \|\widehat{\bu}-{\bu}_1\|_2\max_{1\le j\le p}\|\Theta'_j\|_2=(\lambda_1 v_p+\sigma\sqrt{\log p})\bigg(\frac{\sigma\sqrt{(\lambda_1^2+\sigma^2p)n}}{\lambda_1^2-\lambda_2^2}\land 1 \bigg).
\eeq
To bound the second term, we use H\"older's inequality and the Cauchy-Schwartz inequality
\begin{align*}
\E \max_{1\le j\le p}|(\widehat{\bu}-{\bu}_1)^\top E_j| &\le \E \|\widehat{\bu}-{\bu}_1\|_1\cdot \max_{1\le j\le p}\|E_j\|_\infty\le \sqrt{\E \|\widehat{\bu}-{\bu}_1\|_1^2} \sqrt{\E \max_{\substack{1\le i\le n\\ 1\le j\le p}} |E_{ij}|^2}.
\end{align*}
By (\ref{u.ineq}) and Lemma \ref{pca.lem2}, we have $\E \|\widehat{\bu}-{\bu}_1\|_1^2\le n\E \|\widehat{\bu}-{\bu}_1\|_2^2\le n \big(\frac{\sigma^2{(\lambda_1^2+\sigma^2p)n}}{\lambda_1^2-\lambda_2^2}\land 1 \big)$.
On the other hand, to bound $\E \max_{\substack{1\le i\le n\\ 1\le j\le p}} |E_{ij}|^2$, we need the following lemma concerning a general risk bound for the order statistic of independent Gaussian random variables.

\bep \label{order.risk.lem}
For independent sub-Gaussian random variables $\{y_j\}_{j=1}^p$ with parameter $\sigma^2$ and individual means $\E y_j=\mu_j$, let $y_{(p)}=\max\{y_1,y_2,...,y_p \}$. It holds that, for any $\alpha>0$, $\sup_{\mu_1,...,\mu_p}\E |y_{(p)}-\mu_{(p)}|^\alpha \le C\sigma^\alpha (\log p)^{\alpha/2}$ for some constant $C>0$ only depending on $\alpha$.
\eep

From Proposition \ref{order.risk.lem}, since $E_{ij}$ are independent centered sub-Gaussian random variables with parameter $\sigma^2$, we have $\sqrt{\E \max_{\substack{1\le i\le n, 1\le j\le p}} |E_{ij}|^2} \lesssim \sigma \sqrt{\log p}$.
As a result, we have
\beq \label{homo.max.2}
\E \max_{1\le j\le p}|(\widehat{\bu}-{\bu}_1)^\top E_j|  \lesssim \sigma \sqrt{n\log p}\bigg(\frac{\sigma\sqrt{(\lambda_1^2+\sigma^2p)n}}{\lambda_1^2-\lambda_2^2}\land 1 \bigg).
\eeq
Now for the third term in the last inequality of (\ref{scaler.bnd}), we notice that $\bu_1^\top X_j$ are independent sub-Gaussian with parameter $\sigma^2$ and mean $\lambda v_j$ for $j=1,...,p$. As a result, by Proposition \ref{order.risk.lem}, we have
\beq \label{homo.max.3}
\E \big|\max_{1\le j\le p} \bu_1^\top X_j-\lambda v_p\big|\lesssim \sigma\sqrt{\log p}.
\eeq
Combining (\ref{homo.max.1}) (\ref{homo.max.2}) and (\ref{homo.max.3}), we proved (\ref{homo.max.bnd}). Equation (\ref{homo.upper.eq1}) then follows from (\ref{homo.max.bnd}) (\ref{uhat-u}) and (\ref{homo.upper.decomp}).

\paragraph{Proof of (\ref{homo.upper.eq2}).} As $\frac{1}{p}Y{\bf e}=\frac{1}{p}\Theta {\bf e}+\frac{1}{p}Z {\bf e}={\bf b}+\frac{1}{p}Z{\bf e}$, we have $\E\big\|\frac{1}{p}Y{\bf e}-{\bf b}\big\|_2\le\frac{1}{p} \E\|Z{\bf e}\|_2.$
Since $\frac{1}{p}\|Z{\bf e}\|_2 =\frac{\sigma}{\sqrt{p}} \sqrt{\sum_{i=1}^n\big(\frac{1}{\sigma\sqrt{p}}  \sum_{j=1}^p Z_{ij}\big)^2}$, 
where $\frac{1}{\sigma\sqrt{p}} \sum_{j=1}^p Z_{ij}$ are independent sub-Gaussian random variables.
By the basic property of the sub-Gaussian random variables and the Jensen's inequality, we have $\E \frac{1}{p}\|Z{\bf e}\|_2 \le C\sigma\sqrt{{n}/{p}}$.
This proofs the equation (\ref{homo.upper.eq2}).

\subsection{Minimax Lower Bound: Proof of Theorem \ref{homo.lower.thm}}

The proof of this theorem relies on the following key proposition concerning three important lower bounds over some specially designed parameter spaces.

\bep \label{key.lb.lem}
Suppose $Z$ has i.i.d. entries $Z_{ij}\sim N(0,\sigma^2)$. Then
\begin{itemize}
	\item[(i)] for any $p\ge 8$, $t^2\ge \frac{1-\beta_R^2/2}{\beta_R^2-4/p}\sigma^2\log p$ and $p^{-1/2}\sqrt{\log p}\le \beta_R\le 1$,
	\beq \label{lb.eq1}
	\inf_{\widehat{\Theta}_R}\sup_{ \mathcal{D}_R(t,\beta_R)} \E \|\widehat{\Theta}_R-\Theta_{R}\|_2\gtrsim\sigma\sqrt{\log p};
	\eeq
	\item[(ii)] for any $p\ge 8$, $t^2>0$ and $p^{-1/2}\le \beta_R\le 1$, 
	\begin{align}
	\inf_{\widehat{\Theta}_R}\sup_{\mathcal{D}_{R}(t,\beta_R)} \E \|\widehat{\Theta}_R-\Theta_{R}\|_2&\gtrsim t \beta_R\bigg( 1\land \frac{\sigma\sqrt{ n}}{t}\bigg) \label{lb.eq3};
	\end{align}
	\item[(iii)] for any $c_1p^{-1/2}\sqrt{\log p}\le \beta_R\le c_2<1$ and $\frac{\beta_R^2\sigma^2 p}{1-\beta_R^2}\le t^2\le \sigma^2p/4$, 
	\beq \label{lb.eq4}
	\inf_{\widehat{\Theta}_R}\sup_{\mathcal{D}_{R}(t,\beta_R)} \E \|\widehat{\Theta}_R-\Theta_{R}\|_2\gtrsim t \beta_R\bigg( 1\land \frac{\sigma\sqrt{n(t^2+\sigma^2p)}}{t^2}\bigg);
	\eeq
\end{itemize}
\eep

For any pair $(t,\beta_R)$ satisfying the condition of Theorem \ref{homo.lower.thm}, by Proposition \ref{key.lb.lem} (i), we have $\inf_{\widehat{\Theta}_R}\sup_{\mathcal{D}_{R}(t,\beta_R)} \E \|\widehat{\Theta}_R-\Theta_{R}\|_2\gtrsim \sigma\sqrt{\log p}$, or
\beq \label{psi.rate}
\inf_{\widehat{\Theta}_R}\sup_{\mathcal{D}_{R}(t,\beta_R)} \E \|\widehat{\Theta}_R-\Theta_{R}\|_2\gtrsim \sigma\sqrt{n}\psi(n,p).
\eeq
It remains to prove 
\beq \label{1st.rate}
\inf_{\widehat{\Theta}_R}\sup_{\mathcal{D}_{R}(t,\beta_R)} \E \|\widehat{\Theta}_R-\Theta_{R}\|_2\gtrsim t \beta_R\bigg( 1\land \frac{\sigma\sqrt{n(t^2+\sigma^2p)}}{t^2}\bigg).
\eeq
If $t^2\ge \sigma^2 p/4$, we have ${\sigma\sqrt{n}}/{t} \gtrsim {\sigma\sqrt{n(t^2+\sigma^2p)}}/{t^2},$
so that (\ref{1st.rate}) holds by (\ref{lb.eq3}). If $t^2\le  \sigma^2 p/4$, Proposition \ref{key.lb.lem} (iii) applies directly.  The rest of the proof is devoted to Proposition \ref{key.lb.lem}.

\paragraph{Proof of Proposition \ref{key.lb.lem}.} To establish sharp minimax lower bounds, we consider subsets of $\mathcal{D}$ whose geometric features best describe the essential difficulties for estimating $\Theta_R$.  In particular, we need the following two lemmas concerning testing two and multiple composite hypotheses about an arbitrary operator $F:\Theta\to (\R^d,d)$ characterized by some priors over the parameter space. These general lower bounds can be of independent interest  by providing an effective estimation of  an arbitrary high-dimensional (nonlinear) operator.

\bel[Generalized Le Cam's Method] \label{fuzzy.lem}
Let $\mu_0$ and $\mu_1$ be two priors on the parameter space $\Theta$ of the family $\{P_\theta\}$, and let $P_j$ be the posterior probability measures on $(\mathcal{X},\mathcal{A})$ such that
\[
P_j(S)=\int P_\theta(S)\mu_j(d\theta),\quad \forall S\in\mathcal{A},\quad j=0,1.
\]
Let $F:\Theta\to (\R^d,d)$. If (i) there exist $F_0\in \R^d$ and a $2s$-ball $\mathbb{B}_d(F_0,2s)$ centered at $F_0$ such that $\mu_0(\theta\in\Theta: F(\theta)= F_0) = 1$, $\mu_1(\theta\in\Theta:F(\theta)\notin \mathbb{B}(F_0,2s) )=1,$ and (ii) $\chi^2(P_1,P_0)\le \alpha<\infty$, then
\[
\inf_{\hat{F}}\sup_{\theta\in\Theta} P_{\theta}(d(\hat{F},F(\theta))\ge s)\ge \max\bigg\{ \frac{e^{-\alpha}}{4},\frac{1-\sqrt{\alpha/2}}{2}  \bigg\}.
\]
\eel

\bel[Generalized Fano's Method] \label{fuzzy.lem.2}
Let $\mu_0,\mu_1,...,\mu_M$ be $M+1$ priors on the parameter space $\Theta$ of the family $\{P_\theta\}$, and let $P_j$ be the posterior probability measures on $(\mathcal{X},\mathcal{A})$ such that
\[
P_j(S)=\int P_\theta(S)\mu_j(d\theta),\quad \forall S\in\mathcal{A},\quad j=0,1,...,M.
\]
Let $F:\Theta\to (\R^d,d)$. If (i) there exist some sets $B_0,B_1,...,B_M\subset \R^d$ such that $d(B_i,B_j)\ge 2s$ for some $s>0$ for all $0\le i,j\le M$ and $\mu_j(\theta\in\Theta: F(\theta)\in B_j)= 1$, and (ii) $\frac{1}{M}\sum_{j=1}^MD(P_j,P_0)\le \alpha \log M$ 	with $0<\alpha<1/8$, then
\[
\inf_{\hat{F}}\sup_{\theta\in\Theta} P_{\theta}(d(\hat{F},F(\theta))\ge s)\ge \frac{\sqrt{M}}{1+\sqrt{M}}\bigg(1-2\alpha-\sqrt{\frac{2\alpha}{\log M}} \bigg).
\]
\eel

\paragraph{Proof of (\ref{lb.eq1})}
The proof relies on our general lower bound for testing two composite hypotheses about some high-dimensional operator (Lemma \ref{fuzzy.lem}). We define the one-element set 
\[
\mathcal{H}_0(\eta)=\{ (\Theta,\pi)=({\bf e}_1\bfeta^\top,id) \},
\]
where ${\bf e}_1=(1,0,0,...,0)^\top \in \R^n$ and $\bfeta\in \R^p$ is defined such that, for some $\eta>0$ to be determined later, the first $k=\lfloor p/2\rfloor$ coordinates of $\bfeta$ are $-(\eta,...,\eta)$, the last $k$ coordinates of $\bfeta$ are $(\eta,...,\eta)$, and set the $(\lfloor p/2\rfloor+1)$-th coordinate of $\bfeta$ as 0 if $p$ is odd.
In addition, we define the set
\[
\mathcal{H}_1(\eta)=\{ (\Theta,\pi):  \Theta={\bf e}_1\bfeta'^\top, \pi\in T \},
\]
where $\bfeta'=\bfeta+\frac{\sigma}{2}\sqrt{\log p}\cdot (-1,0,...,0,1)^\top$, and $T=\{\pi\in \mathcal{S}_p:$ $\pi$ is a transposition between $k$-th and $p$-th element, $\forall \lfloor p/2\rfloor+1 \le k\le p \}.$
For any pair $(t,\beta_R)$ satisfying the condition of Proposition \ref{key.lb.lem} (i), to find specific $\eta$ such that $\mathcal{H}_0(\eta)\subset \mathcal{D}_{R}(t,\beta_R)$ and $\mathcal{H}_1(\eta)\subset \mathcal{D}_{R}(t,\beta_R)$, we consider $\eta$ such that
\beq \label{check.cond.1}
t=\|{\bf e}_1\|_2\|\bfeta\|_2=\sqrt{p\eta^2},\quad \beta_R^2\ge\frac{2(\eta+\frac{\sigma}{2}\sqrt{\log p})^2}{\|\bfeta'\|_2^2}\ge \frac{\eta^2}{\|\bfeta\|_2^2},
\eeq
where the last inequality holds by construction. Note that the second equation holds if
\[
\beta_R^2\ge\frac{4\eta^2+\sigma^2\log p}{p\eta^2+\frac{\sigma^2}{2}\log p}\quad\text{or}\quad
\eta^2 \ge\frac{1-\beta_R^2/2}{p\beta^2_R-4}\cdot \sigma^2\log p.
\]
Plug this into the first equation of (\ref{check.cond.1}), we have $t^2=p\eta^2\ge\frac{1-\beta_R^2/2}{\beta^2_R-4p^{-1}}\cdot \sigma^2\log p$.
Hence, as long as $t^2\ge \frac{1-\beta_R^2/2}{\beta_R^2-4/p}\sigma^2\log p$, we can always find $\eta$ such that $\mathcal{H}_0(\eta)$. On the other hand, it can be checked that for any $(\Theta,\pi)\in \mathcal{H}_1(\eta)$,
\[
t^2\le p\eta^2+\frac{\sigma^2}{2}\log p\le\|\Theta\|^2\le (p+2)\eta^2+\sigma^2\log p\le 4t^2,
\]
so that $\mathcal{H}_1(\eta)\subset \mathcal{D}_{R}(t,\beta_R)$.

Intuitively,  the construction of $\mathcal{H}_0$ and $\mathcal{H}_1$ reflects the effects of perturbing the extreme values of the first right singular vector of $\Theta$. The subset $\mathcal{H}_1$ is constructed so that a mixture distribution based on parameters uniformly chosen from $\mathcal{H}_1$ will be statistically indistinguishable from the null distribution defined on $\mathcal{H}_0$. 
Toward this end, we define $\mu_1$ as the uniform prior over $\mathcal{H}_1(\eta)$, and define $\mu_0$ as the point-mass at the one-point set $\mathcal{H}_0(\eta)$. Let $F: \mathcal{D}\to \R^n$ be the function such that $F(\Theta)$ is the rightmost column $\Theta_R$ of $\Theta$. It holds that
\[
\mu_0((\Theta,\mu)\in \mathcal{D}: F(\Theta)=(\eta,0,...,0)^\top \in\R^n  )=1,
\]
\[
\mu_1((\Theta,\pi)\in \mathcal{D}: F(\Theta)= (\eta+\frac{\sigma}{2}\sqrt{\log p},0,...,0)^\top \in \R^n )=1.
\]
Then $\| F(\Theta_0)-F(\Theta_1)\|_2=\frac{\sigma}{2}\sqrt{\log p}$ for any $\Theta_0\in \mathcal{H}_0$ and $\Theta_1\in \mathcal{H}_1$. Moreover, let $P_0$ and $P_1$ be the posterior distributions of $Y=\Theta\Pi+Z$. The following lemma controls the chi-square divergence $\chi^2(P_1,P_0)$.

\bel \label{chisq.lem}
Under the condition of Proposition \ref{key.lb.lem} (i), it holds that $\chi^2(P_1,P_0) \le 3$.
\eel

Lemmas \ref{fuzzy.lem} and \ref{chisq.lem} yield $\inf_{\widehat{\Theta}}\sup_{\mathcal{D}_R(t,\beta_R)} P\big( \|\widehat{\Theta}-\Theta_{R}\|_2 \ge \frac{\sigma}{4}\sqrt{\log p}\big) \ge 0.01$.
The final result then follows from the Markov's inequality.

\paragraph{Proof of (\ref{lb.eq3})}  The proof of this lower bound relies on the following general lower bound for testing multiple hypotheses.

\bel[\cite{tsybakov2009introduction}] \label{lower.lem}
Assume that $M\ge 2$ and suppose that $\Theta$ contains elements $\theta_0,\theta_1,...,\theta_M$ such that: (i) $d(\theta_j,\theta_k)\ge 2s>0$ for any $0\le  j<k\le M$; (ii) it holds that $\frac{1}{M}\sum_{j=1}^M D(P_j,P_0)\le \alpha \log M$ with $0<\alpha<1/8$ and $P_j=P_{\theta_j}$ for $j=0,1,...,M$. Then
\[
\inf_{\hat{\theta}}\sup_{\theta\in\Theta}P_\theta(d(\hat{\theta},\theta)\ge s) \ge \frac{\sqrt{M}}{1+\sqrt{M}}\bigg( 1-2\alpha-\sqrt{\frac{2\alpha}{\log M}}\bigg) >0.
\]
\eel 

We consider some subset $\mathcal{D}_m\subset \mathcal{D}(t,\beta_R)$ for any $(t,\beta_R)$ satisfying the condition of Proposition \ref{key.lb.lem} (ii). Specifically, for some fixed unit vector $\bv_0=(v_1,...,v_p)^\top \in \R^p$ where $v_1\le v_2\le ...\le v_p=\beta_R$ and $\sum_{j=1}^p v_j=0$, we define
\[
\mathcal{D}_{\bv_0}=\{ (\Theta_\bu,id)\in \mathcal{D}_R(t,\beta_R): \Theta_\bu=t\bu\bv_0^\top, \|\bu\|_2=1\}.
\]
For any $(\Theta_\bu,id)\in \mathcal{D}_\bv$, let $P_\bu$ be the probability measure of $Y=\Theta_\bu+Z\in \R^{n\times p}$ with $Z_{ij}\sim N(0,\sigma^2)$. Then the KL-divergence satisfies $D(P_\bu\| P_{\bu'}) = \frac{1}{2\sigma^2}\|t\bu\bv_0^\top-t\bu'{\bv_0}^{\top}\|_F^2=\frac{t^2}{2\sigma^2}\|\bu-\bu'\|_2^2$.
Now for any fixed unit vector $\bu_0\in \R^n$, we consider the ball with radius $0<\epsilon<1$ and centered at $\bu_0$
\[
B(\bu_0,\epsilon)=\{\bu\in \R^n: \|\bu\|_2=1, \sqrt{1-(\bu_0^\top \bu)^2} \le \epsilon \}
\]
Since $\sqrt{1-(\bu^\top \bu')^2}=\|\bu\bu^\top-{\bu'}{\bu'}^{\top}\|_F/\sqrt{2}$, to construct a local packing of $B(\bu_0,\epsilon)$, we can use the following lemma regarding the metric entropy of the  Grassmannian manifold $G(k,r)$ from \cite{cai2013sparse} (see also Proposition 8 of \cite{szarek1982nets}.) 

\bel[\cite{cai2013sparse}] \label{entropy.g}
For any $V\in O(k,r)$, identifying the subspace $\text{span}(V)$ with its projection matrix $VV^\top$, define the metric on the Grassmannian manifold $G(k,r)$ by $\rho(VV^\top, UU^\top)=\|VV^\top-UU^\top\|_F$. Then for any $\epsilon\in (0,\sqrt{2(r\land (k-r))})$, we have $\big( \frac{c_0}{\epsilon}\big)^{r(k-r)} \le \mathcal{N}(G(k,r),\epsilon)\le \big( \frac{c_1}{\epsilon}\big)^{r(k-r)}$,
where $\mathcal{N}(E,\epsilon)$ is the $\epsilon$-covering number of $E$ and $c_0,c_1$ are absolute constants. Moreover, for any $V\in O(k,r)$ and any $\alpha\in(0,1)$, $\mathcal{M}(B(V,\epsilon),\alpha\epsilon) \ge \big(\frac{c_0}{\alpha c_1} \big)^{r(k-r)}$,
where $\mathcal{M}(E,\epsilon)$ is the $\epsilon$-packing number of $E$.
\eel

As a result, for any $\alpha,\epsilon \in (0,1)$, there exist $\bu_1,...,\bu_m\in B(\bu_0,\epsilon)$, such that $m\ge (c/\alpha)^{n-1}$ and $\min_{1\le i\ne j\le m}\sqrt{1-(\bu_i^\top \bu_j)^2}\ge \alpha\epsilon$. 
Without loss of generality, we choose $\bu_1,...,\bu_m$ such that $\bu_i^\top \bu_j\ge 0$ for all $i,j\in \{1,...,m\}$ (this can be done if $\epsilon^2<0.5$). Hence, by focusing on the $m$ point subset $\mathcal{D}_m=\{ (\Theta_\bu,id)\in \mathcal{D}_{\bv_0}: \bu\in \{ \bu_1,...,\bu_m \}\}$, for any $1\le i\ne j\le m$,
\begin{align*}
&\|F(\Theta_{\bu_i})-F(\Theta_{\bu_j})\|_2=
\|t\bu_i\beta_R-t\bu_j\beta_R\|_2 \ge t\beta_R\sqrt{2(1-\bu_i^\top \bu_j)}\\
&=t\beta_R\sqrt{\frac{2(1-(\bu_i^\top \bu_j)^2)}{(1+\bu_i^\top \bu_j)}}\ge t\beta_R\sqrt{1-(\bu_i^\top \bu_j)^2}\ge t\beta_R\alpha\epsilon
\end{align*}
and the KL-divergence $\max_{1\le i\le m} D(P_{\bu_i}\|P_{\bu_0}) = \frac{t^2}{2\sigma^2}\|\bu_i-\bu_0\|_2^2\le  t^2 \epsilon^2/\sigma^2.$
If we set $\epsilon^2=c(\sigma^2n/t^2 \land 1)$ for some constant $c>0$, we have $\max_{1\le i\le m} D(P_{\bu_i}\|P_{\bu_0}) \le \frac{ct^2}{\sigma^2}\big( \frac{\sigma^2 n}{t^2}\land 1 \big)\le \log m.$
Then Lemma \ref{lower.lem} yields $\inf_{\widehat{\Theta}}\sup_{\mathcal{D}_m} P\big( \|\widehat{\Theta}-\Theta_{R}\|_2\ge ct\beta_R\epsilon \big)\ge C$ for some constants $C,c>0$. The desired lower bound follows from the Markov inequality and the inclusion $\mathcal{D}_m\subset \mathcal{D}_{\bv_0}\subset \mathcal{D}_{R}(t,\beta_R)$.

\paragraph{Proof of (\ref{lb.eq4})} The proof of the last lower bound is more complicated than the previous ones. Specifically, we need to construct mixture distributions that reflect the uncertainty of both the first right and the first left singular vectors of $\Theta$. Throughout, we set $\Pi={\bf I}$.

Recall that $\beta_R \in [c_1p^{-1/2}\sqrt{\log p},1]$. We define the following class of density $P_Y$ of $Y=\Theta+Z$, where $Z_{ij}\sim_{i.i.d.} N(0,\sigma^2)$ and $\Theta=\lambda \bu \bv^\top$ for some fixed unit vectors $\bu\in \R^n$ and $\bv\in \R^p$ with some constraints, i.e.,
\[
\mathcal{P}_{\bu,t,\beta_R}=\Bigg\{ P_Y:  \begin{aligned}
&Y=\Theta+Z,  \Theta = \lambda \bu\bv^\top, \frac{\beta_R}{36}\le\frac{v_{(p)}-\bar{v}}{V^{1/2}(\bv)} \le \beta_R, \\
& t\le \lambda\le 6t, \|\bu\|_2=\|\bv\|_2=1 
\end{aligned} \Bigg\}.
\]
where for $\xx=(x_1,...,x_p)$, we denote $\bar{x}=p^{-1}\sum_{j=1}^p x_j$, $V(\xx)=\sum (x_i-\bar{x})^2$ and $x_{(p)}$ is the largest component of $\xx$.
In particular, the permutation matrix $\Pi$ amounts to permuting the relative orders of the coordinates of $\bv$, which is a monotone vector before permutation.
In addition, for any $\bw\in \R^{p-1}$, define $\bw_{+}=\begin{bmatrix} \bw\\ \frac{\beta_R}{4} \end{bmatrix}\in \R^p$,
and define $\mathcal{U}=\big\{ \bu\in \R^n: |u_j|=1/\sqrt{n}    \big\}$.
Denote
\[
\mathcal{G}=\{ \bw\in \R^{p-1}: 1/2\le \|\bw\|_2\le 2,\\\max\{w_{(p)}, |\bar{w}|\}\le \frac{c_1}{8}p^{-1/2}\sqrt{\log p} \},
\]
where $w_{(p)}$ is the largest component of $\bw$.
We construct the following Gaussian mixture measure 
\begin{align*}
\bar{P}_{\bu,t,\beta_R}(Y)=C_{t,\beta_R}\int_{\mathcal{G}} & \frac{\sigma^{np}}{(2\pi)^{np/2}}\exp(-\|Y-2t\bu\bw_{+}^\top\|_F^2/(2\sigma^2))\\
&\times \bigg( \frac{p-1}{2\pi}\bigg)^{(p-1)/2}\exp(-(p-1)\|\bw\|_2^2/2)dw
\end{align*}
Here $C_{t,\beta_R}$ is the constant which normalizes the integral and makes $\bar{P}_{\bu,t,\beta_R}$ a valid probability density. To be specific, $C^{-1}_{t,\beta_R} = P\big( \bw\in \mathcal{G} | w_j\sim N(0,(p-1)^{-1}), 1\le j\le p-1  \big)$.
Moreover, since  in the event $\mathcal{G}$, $2t\bu\bw_{+}^\top$ is rank one with the singular value $2t\|\bw_{+}\|_2\in (t,2\sqrt{5}t)$ and by construction we have $\frac{\beta_R}{36}\le\frac{w_{+,(p)}-\bar{w}_{+}}{V^{1/2}(\bw_+)}\le  \beta_R$,
it follows that
$\bar{P}_{\bu,t,\beta_R}(Y)$ is a mixture density of infinite members of $\mathcal{P}_{\bu,t,\beta_R}$, i.e., $\bar{P}_{\bu,t,\beta_R}(Y)\in \text{Conv}(\mathcal{P}_{\bu,t,\beta_R})$.

The rest of the proofs rely on the our general lower bound based on testing multiple composite hypotheses (Lemma \ref{fuzzy.lem.2}), as well as the following lemma that gives an upper bound for the KL-divergence between $\bar{P}_{\bu,t,\beta_R}$ and $\bar{P}_{\bu',t,\beta_R}$ for any $\bu,\bu'\in \R^n$.

\bel \label{kl.mix.prop}
Under the assumption of Proposition \ref{key.lb.lem} (iii), for any unit vectors $u,u'\in \R^n$, we have $D(\bar{P}_{\bu,t,\beta_R}\| \bar{P}_{\bu',t,\beta_R}) \le \frac{C_1t^4}{\sigma^2(4t^2+\sigma^2(p-1))}(1-(\bu^\top \bu')^2)+C_{2}$,
where $C_{1},C_2>0$ are some uniform constant.
\eel

Again by Lemma \ref{entropy.g}, for any $\alpha,\epsilon \in (0,1)$, there exists $\bu_1,...,\bu_m\in B(\bu_0,\epsilon)$, such that $m\ge (c/\alpha)^{n-1}$ and $\min_{1\le i\ne j\le m}\sqrt{1-(\bu_i^\top \bu_j)^2}\ge \alpha\epsilon$. 
Again, we choose $\bu_1,...,\bu_m$ such that $\bu_i^\top \bu_j\ge 0$ for all $i,j\in \{1,...,m\}$. Hence, for $i=1,...,m$, let $\mu_i$ be the priors over the parameter space of $\Theta$ leading to the posterior $\bar{P}_{\bu_i,t,\beta_R}$. It holds that, for any $\Theta_i\in \text{supp}(\mu_i)$ and $\Theta_j\in \text{supp}(\mu_j)$ with $1\le i\ne j\le m$, $\|F(\Theta_i)-F(\Theta_j)\|_2\ge C
\|tu_i\beta_R-tu_j\beta_R\|_2 \ge Ct\beta_R\sqrt{2(1-u_i^\top u_j)}\ge Ct\beta_R\alpha\epsilon$,
and the KL-divergence $\max_{1\le i\le m} D(P_{\bu_i,t,\beta_R}\|P_{\bu_0,t,\beta_R}) \le \frac{C_1 t^4\epsilon^2}{\sigma^2(4t^2+\sigma^2(p-1))}+C_2.$
Now set $\epsilon^2= c\big(\frac{\sigma^2n(t^2+\sigma^2 p)}{t^4}\land 1\big)$ for some constant $c>0$. Lemma \ref{fuzzy.lem.2} and the Markov inequality yield
\beq
\inf_{\widehat{\Theta}} \sup_{\supp(\mu_0),...,\supp(\mu_m)} \E \|\widehat{\Theta}-\Theta_R\|_2\gtrsim t\beta_R\bigg(  \frac{\sigma\sqrt{(t^2+\sigma^2 p)n}}{t^2}\land 1 \bigg),
\eeq
where the supremum is over the union of all the sets $\{\text{supp}(\mu_i)\}_{i=0}^m$.
Finally, since $\cup_{i=0}^m\supp(\mu_i)\times \{id\}\subset \mathcal{D}_{R}(t,\beta_R)$, equation (\ref{lb.eq4}) holds.

\section{Optimal Range Vector Estimation: Proof of Theorem \ref{r.homo.upper.thm}}

\paragraph{Risk upper bounds.} Observe that
\begin{align*}
\|\widehat{\bf R}^*-{\bf R} \|_2&=\| (\widehat{\Theta}^*_R-\widehat{\Theta}_L^*)-(\Theta_R-\Theta_L)  \|_2\le \|\widehat{\Theta}_{R}^*-\Theta_{R}\|_2+\|\widehat{\Theta}_{L}^*-\Theta_{L}\|_2.
\end{align*}
The upper bound then follows by applying Theorem \ref{unif.general.thm} for both $\Theta_R$ and $\Theta_L$.

\paragraph{Minimax lower bounds.} Similar to the proof of Theorem \ref{homo.lower.thm}, we apply the following proposition which can be obtained by directly generalizing the proofs of Proposition \ref{key.lb.lem} of our main paper.

\bep 
Suppose $Z$ has i.i.d. entries $Z_{ij}\sim N(0,\sigma^2)$. Then
\begin{itemize}
	\item[(i)] for any $p\ge 8$, $t^2\ge \frac{1-(\beta_R^2\land \beta_L^2)/2}{(\beta_R^2\land \beta_L^2)-4/p}\sigma^2\log p$ and $p^{-1/2}\sqrt{\log p}\le \beta_R,\beta_L\le 1$, 
	\beq \label{lb2.eq1}
	\inf_{\widehat{\bf R}}\sup_{ \mathcal{D}_W(t,\beta_R,\beta_L)} \E \|\widehat{\bf R}-{\bf R}\|_2\gtrsim \sigma\sqrt{\log p};
	\eeq
	\item[(ii)] for any $p\ge 8$, $t^2> 0$ and $p^{-1/2}\le \beta_R,\beta_L\le 1$, 
	\begin{align}
	\inf_{\widehat{\bf R}}\sup_{\mathcal{D}_{W}(t,\beta_R,\beta_L)} \E \|\widehat{\bf R}-{\bf R}\|_2&\gtrsim t (\beta_R+\beta_L)\bigg( 1\land \frac{\sigma\sqrt{ n}}{t}\bigg) \label{lb2.eq3};
	\end{align}
	\item[(iii)] for any $\frac{(\beta_R^2\lor \beta_L^2)\sigma^2 p}{1-(\beta_R^2\lor \beta_L^2)}\le t^2\le \sigma^2p/4$ and $c_1p^{-1/2}\sqrt{\log p}\le\beta_R,\beta_L\le c_2<1$, 
	\beq \label{lb2.eq4}
	\inf_{\widehat{\bf R}}\sup_{\mathcal{D}_{W}(t,\beta_R,\beta_L)} \E \|\widehat{\bf R}+{\bf R}\|_2\gtrsim t (\beta_R-\beta_L)\bigg( 1\land \frac{\sqrt{n(t^2+\sigma^2p)}}{t^2}\bigg);
	\eeq
\end{itemize}
\eep

With these lower bounds, the proof of the minimax lower bounds in Theorem \ref{r.homo.upper.thm} follows the same argument as the proof of Theorem \ref{homo.lower.thm}.

\section{Equivalent Estimators: Proof of Equation (\ref{equiv}) of the Main Paper}

Let $\hat{\Pi}$ be the permutation matrix corresponding to the permutation $\hat{\pi}$.
For $i=1,...,n$, we have $\check{Y}_{i.}=(y_{i,\hat{\pi}(1)},...,y_{i,\hat{\pi}(p)})=Y_{i.}\hat{\Pi}^{-1}$. Similarly, $(\widehat{v}_{\hat{\pi}(1)},\widehat{v}_{\hat{\pi}(2)},...,\widehat{v}_{\hat{\pi}(p)})=\widehat{\bv}^\top\hat{\Pi}^{-1}$. Fitting a simple linear regression between $\check{Y}_{i.}$ and $(\widehat{v}_{\hat{\pi}(1)},\widehat{v}_{\hat{\pi}(2)},...,\widehat{v}_{\hat{\pi}(p)})$ is therefore the same as fitting a regression between $Y_{i.}$ and $\widehat{\bv}^\top$. As a result, if we denote $m(\widehat{\bv})=\frac{1}{p}\sum_{j=1}^p \widehat{v}_j$, then
\[
\beta_i=\frac{\sum_{j=1}^p (\widehat{v}_j-m(\widehat{\bv}))(y_{ij}-\bar{y}_i)}{\sum_{j=1}^p (\widehat{v}_j-m(\widehat{\bv}))^2}=\frac{\sum_{j=1}^p\widehat{v}_j(y_{ij}-\bar{y}_i)}{\sum_{j=1}^p (\widehat{v}_j-m(\widehat{\bv}))^2}=\frac{\widehat{v}^\top X_{i.}}{\sum_{j=1}^p (\widehat{v}_j-m(\widehat{\bv}))^2},
\]
and $\alpha_i=\bar{Y}_{i.}-m(\widehat{\bv})\beta_i=\frac{1}{p}{\bf e}^\top Y_{i.}-m(\widehat{\bv})\beta_i.$
Hence
\beq \label{ThetaRReg}
\hat{\Theta}_{R}^{Reg}=\frac{(\widehat{v}_{(p)}-m(\widehat{\bv})X\widehat{\bv}}{\sum_{j=1}^p (\widehat{v}_j-m(\widehat{\bv}))^2}+\frac{1}{p}Ye=\frac{(\widehat{v}_{(p)}-m(\widehat{\bv})X\widehat{\bv}}{1-pm^2(\widehat{\bv})}+\frac{1}{p}Y{\bf e}.
\eeq
Now recall that $X=Y({\bf I}_p-\frac{1}{p}{\bf ee}^\top)$ so that $X {\bf e}=0$. By SVD of $X$, each right singular vector of $X$ is orthogonal to ${\bf e}$, so that $\widehat{\bv}^\top {\bf e}=pm(\widehat{\bv})=0$.  Plugging $m(\widehat{\bv})=0$ to (\ref{ThetaRReg}), we have $\hat{\Theta}_{R}^{Reg}=X\widehat{\bv}+\frac{1}{p}Y{\bf e}=\widehat{\Theta}_R^*$.
Similar arguments can be applied to prove $	\hat{\Theta}_{L}^{Reg}=\widehat{\Theta}_L^*$ and $\hat{\bf R}^{Reg}=\widehat{\bf R}^*$.

\section{Proofs of Technical Results}

\subsection{Properties of the Monotone Matrices: Proof of Proposition \ref{rsv.prop}}

We first show the second statement that the vector $\text{sgn}(\bu_1)$ indicates the direction of monotonicity of the rows of $\Theta'$.

By SVD of $\Theta'=(\theta'_{ij})\in \R^{n\times p}$, its first left singular vector
\beq 
\bu_1 = \argmax_{\|x\|_2=1} x^\top \Theta'{\Theta'}^\top x=\argmax_{\|x\|_2=1} \sum_{j=1}^p \bigg( \sum_{i=1}^n x_i \theta'_{ij}\bigg)^2.
\eeq
Now let $m(\Theta')=(m_i)_{1\le i\le n} \in \{ -1,1 \}^n$ be the vector of $n$ signed indicators of whether each row of $\Theta'$ is nondecreasing. 
In the following, we show that, for any unit vector $x\in \R^n$,
\beq \label{ineq.2}
\sum_{j=1}^p \bigg( \sum_{i=1}^n x_i \theta'_{ij}\bigg)^2\le \sum_{j=1}^p \bigg( \sum_{i=1}^n m_i|x_i|\cdot \theta'_{ij}\bigg)^2,
\eeq
which combined with the definition of $\bu_1$ immediately implies the second statement. Toward this end, note that
\[
\sum_{j=1}^p \bigg( \sum_{i=1}^n x_i \theta'_{ij}\bigg)^2=\sum_{j=1}^p\sum_{i=1}^nx^2_i{\theta'}^2_{ij}+\sum_{i\ne k}x_ix_k\bigg(\sum_{j=1}^p{\theta'}_{ij}{\theta'}_{kj}\bigg).
\]
Then the inequality (\ref{ineq.2}) follows from
\[
\sum_{i\ne k}x_ix_k\bigg(\sum_{j=1}^p{\theta'}_{ij}{\theta'}_{kj}\bigg)\le \sum_{i\ne k}m_im_k|x_ix_k|\bigg(\sum_{j=1}^p{\theta'}_{ij}{\theta'}_{kj}\bigg),
\]
which is true as long as $m_im_k\sum_{j=1}^p{\theta'}_{ij}{\theta'}_{kj}\ge 0$ for all pairs $i\ne k$. Hence, we conclude the proof of the second statement by showing the following lemma.

\bel
For any nondecreasing vectors $a,b\in \R^n$, such that $\sum_{i=1}^n a_i=0$. Then it follows that $a^\top b\ge 0$.
\eel

\begin{proof}
	Since $a$ and $b$ are both nondecreasing, there exist a constant $\delta$ such that the components of $b+\delta\cdot \bf{1}$ has the same sign as $a$. Hence the claim follows from $0\le a^\top(b+\delta\cdot {\bf{1}}) = a^\top b+\delta a^\top {\bf{1}}=a^\top b.$ (Thank Rui Duan for providing this elegant proof.)
\end{proof}

Next we show the first statement in the proposition.
Note that by SVD of $\Theta'$, we have 
\beq \label{v1.def}
\bv_1 = \argmax_{\|x\|_2=1} x^\top \Theta'^\top \Theta' x=\argmax_{\|x\|_2=1} \sum_{i=1}^n \bigg( \sum_{j=1}^p x_j \theta'_{ij}\bigg)^2.
\eeq
To prove that $\bv_1$ is monotone, we need the following rearrangement inequality.
\bel [Rearrangement Inequality] \label{re.ineq}
If $a_1\ge a_2\ge ...\ge a_n$ and $b_1\ge b_2\ge ... \ge b_n$, then
\[
a_nb_1+...+a_1b_n \le a_{\sigma(1)}b_1+...+a_{\sigma(n)}b_n \le a_1b_1+...+a_nb_n,
\]
where $\sigma$ is any permutation in $S_n$.
\eel

Without loss of generality, we assume $\bu_1=m(\Theta')$, indicating whether each row of $\Theta'$ is nondecreasing.
Now suppose $\bv_1$ is not a monotone vector. Then for any $i=1,...,n$, if $\sum_{j=1}^pv_j\theta'_{ij}\ge 0$, since $\lambda_1(\Theta')u_i=\sum_{j=1}^pv_j\theta'_{ij}$, it means $u_i\ge 0$, or, the $i$th row of $\Theta'$ is nondecreasing. Thus by Lemma \ref{re.ineq}, rearranging the components of $\bv_1$ in nondecreasing order will further increase the value $\big( \sum_{j=1}^p x_j \theta'_{ij}\big)^2$. Similarly, if $\sum_{j=1}^pv_j\theta'_{ij}< 0$, then the $i$th row of $\Theta'$ is decreasing, and thus, by Lemma \ref{re.ineq}, rearranging the components of $\bv_1$ in nondecreasing order will also increase the value $\big( \sum_{j=1}^p x_j \theta'_{ij}\big)^2$. Hence, by definition of $\bv_1$, the components of $\bv_1$ should be of nondecreasing order.

\subsection{Two General Minimax Lower Bounds}

In this section, we prove two general lower bounds for testing composite (fuzzy) hypotheses, namely Lemmas \ref{fuzzy.lem} and \ref{fuzzy.lem.2}, which generalize the Le Cam's method and the Fano's method to an arbitrary high-dimensional operator.

\subsubsection{Proof of Lemma  \ref{fuzzy.lem} }
Observe that
\begin{align*}
\int P_\theta(d(\hat{F},F(\theta))\ge s) \mu_0(d\theta)&\ge \int I(d(\hat{F},F_0)\ge s, F(\theta)=F_0) dP_\theta\mu_0(d\theta)\\
&\ge \int I(d(\hat{F},F_0)\ge s) dP_\theta\mu_0(d\theta)-\int I(F(\theta)\ne F_0)dP_\theta\mu_0(d\theta)\\
&=P_0(d(\hat{F},F_0)\ge s).
\end{align*}
In a similar way,
\begin{align*}
\int P_\theta(d(\hat{F},F(\theta))\ge s) \mu_1(d\theta)&\ge \int I(\hat{F} \in  \mathbb{B}(F_0,s), F(\theta)\notin  \mathbb{B}(F_0,2s) ) dP_\theta\mu_1(d\theta)\\
&\ge \int I(\hat{F} \in  \mathbb{B}(F_0,s)) dP_\theta\mu_1(d\theta)-\int I(F(\theta)\in  \mathbb{B}(F_0,2s))dP_\theta\mu_1(d\theta)\\
&\ge P_1(d(\hat{F},F_0)< s).
\end{align*}
Now since
\begin{align*}
\sup_{\theta\in\Theta} P_\theta( d(\hat{F},F(\theta))\ge s) \ge \max\bigg\{ \int P_\theta( d(\hat{F},F(\theta))\ge s )\mu_0(d\theta), \int P_\theta( d(\hat{F},F(\theta))\ge s )\mu_1(d\theta) \bigg\},
\end{align*}
we have
\begin{align*}
\sup_{\theta\in \Theta} P_\theta(d(\hat{F},F(\theta))\ge s )& \ge \max\bigg\{  P_0(d(\hat{F},F_0)\ge s), P_1(d(\hat{F},F_0)< s) \bigg\}\\
&= \inf_{\psi}\max\{ P_0(\psi =1),P_1(\psi =0)  \}\\
&\ge \inf_{\psi} \frac{P_0(\psi =1)+P_1(\psi =0)  }{2}\\
&=\frac{1}{2}\int \min(dP_0,dP_1 )
\end{align*}
where $\inf_{\psi}$ denotes the infimum over all tests $\psi$ taking values in $\{0,1\}$. The proof is complete by using Theorem 2.2 of \cite{tsybakov2009introduction}.

\subsubsection{Proof of Lemma  \ref{fuzzy.lem.2}}

Observe that, for any $j=0,...,M$,
\begin{align*}
\int P_\theta(d(\hat{F},F(\theta))\ge s) \mu_j(d\theta)&\ge \int I(d(\hat{F},B_j)\ge s, F(\theta)\in B_j) dP_\theta\mu_j(d\theta)\\
&\ge \int I(d(\hat{F},B_j)\ge s) dP_\theta\mu_j(d\theta)\\
&\ge P_j(d(\hat{F},B_j)\ge s).
\end{align*}
Now since
\begin{align*}
\sup_{\theta\in\Theta} P_\theta( d(\hat{F},F(\theta))\ge s) \ge \max\bigg\{ \int P_\theta( d(\hat{F},F(\theta))\ge s )\mu_0(d\theta), \max_{1\le j\le M}\int P_\theta( d(\hat{F},F(\theta))\ge s )\mu_j(d\theta) \bigg\},
\end{align*}
we have
\begin{align}
\sup_{\theta\in \Theta} P_\theta(d(\hat{F},F(\theta))\ge s )& \ge \max\bigg\{  P_0(d(\hat{F},B_0)\ge s)-\beta_0, \max_{1\le j\le p}[P_j(d(\hat{F},B_j)>s)-\beta_j] \bigg\} \nonumber \\
&\ge \inf_{\psi}\max\{ P_0(\psi \ne 0),\max_{1\le j\le p}[P_j(\psi \ne j) ] \} \nonumber \\
&\ge  \inf_{\psi}\max\bigg\{ P_0(\psi \ne 0),\frac{1}{M}\sum_{j=1}^pP_j(\psi \ne j) \bigg\}, \label{m.lb.eq}
\end{align}
where  $\inf_{\psi}$ denotes the infimum over all tests $\psi$ taking values in $\{0,1,...,M\}$. Let $P_{0,j}^a$ be the absolutely continuous component of the measure $P_0$ with respect to $P_j$, and let the random even $A_j=\bigg\{ \frac{dP_{0,j}^a}{dP_j}\ge \tau \bigg\}$  for $j=1,...,M$. We write
\begin{align}
P_0(\psi\ne 0)&=\sum_{j=1}^M P_0(\psi=j)\ge \sum_{j=1}^M P_{0,j}^a(\psi=j)\nonumber \\
&\ge  \sum_{j=1}^M \tau P_{j}(\{\psi=j \}\cap A_j) \nonumber \\
&\ge \tau M\bigg( \frac{1}{M}\sum_{j=1}^M P_j(\psi=j) \bigg)-\tau \sum_{j=1}^MP_j(A_j^c). \label{p0.eq}
\end{align}
Combining (\ref{m.lb.eq}) and (\ref{p0.eq}), we have
\begin{align*}
&\sup_{\theta\in \Theta} P_\theta(d(\hat{F},F(\theta))\ge s )\\
&\ge  \inf_{\psi}\max\bigg\{ M\bigg( \frac{1}{M}\sum_{j=1}^M P_j(\psi=j) \bigg)-\tau \sum_{j=1}^MP_j(A_j^c) ,\frac{1}{M}\sum_{j=1}^pP_j(\psi \ne j) \bigg\}\\
&\ge \min_{0\le r\le 1}\max\bigg\{ \tau M\bigg(r-\frac{1}{M} \sum_{j=1}^MP_j(A_j^c)\bigg) ,1-r \bigg\}\\
&=\frac{\tau M}{\tau M+1}\bigg( \frac{1}{M} \sum_{j=1}^MP_j\bigg(\frac{dP_{0,j}^a}{dP_j}\ge \tau \bigg) \bigg).
\end{align*}
In the following, we check that for all $0<\tau<1$,
\beq \label{llr}
\frac{1}{M}\sum_{j=1}^MP_j \bigg(\frac{dP_{0,j}^a}{dP_j}\ge \tau \bigg) \ge 1+\frac{\alpha \log M+\sqrt{\alpha/2\cdot\log M}}{\log \tau}.
\eeq
Recall that 
\beq \label{KL}
\frac{1}{M}\sum_{j=1}^MD(P_j,P_0)\le \alpha \log M,
\eeq
we have $P_j\ll P_0$ and $dP_j/dP_0=dP_j/dP_{0,j}^a$ everywhere except for a set having $P_j$-measure zero. Then we have
\begin{align*}
P_j \bigg(\frac{dP_{0,j}^a}{dP_j}\ge \tau \bigg) &=P_j \bigg(\frac{dP_j}{dP_0}\le \frac{1}{\tau} \bigg) \\
&=1-P_j \bigg(\log \frac{dP_j}{dP_0}>\log \frac{1}{\tau}  \bigg) \\
&\ge 1-\frac{1}{\log(1/\tau)} \int \bigg( \log \frac{dP_j}{dP_0}\bigg)_{+}dP_j\quad [\text{Markov's inequality}]\\
&\ge 1- \frac{1}{\log(1/\tau)} \big[ D(P_j,P_0)+\sqrt{D(P_j,P_0)/2} \big] \quad [\text{2nd Pinsker's inequality}].
\end{align*}
Note that the second Pinsker's inequality can be found in \cite{tsybakov2009introduction}.
Now by the Jensen's inequality and (\ref{KL}), we have
\[
\frac{1}{M}\sum_{j=1}^M \sqrt{D(P_j,P_0)} \le \bigg(\frac{1}{M}\sum_{j=1}^M D(P_j,P_0) \bigg)^{1/2}\le \sqrt{\alpha \log M},
\]
which proves (\ref{llr}). Finally, the proof is completed by taking $\tau=M^{-1/2}$.

\subsection{An Elementary Inequality: Proof of Lemma \ref{max.ineq.lem}}

The proof is separated into three cases.

\emph{Case 1}: Suppose both $\max_i a_i$ and $\max_i b_i$ are non-negative. In this case, without loss of generality, we assume $\max_i a_i\ge \max_i b_i\ge 0$. Let $1\le k_0\le n$ such that $a_{k_0}=\max_i a_i$. Then apparently
\[
a_{k_0}-b_{k_0}=|a_{k_0}-b_{k_0}|\le \max_{1\le i\le n}|a_i-b_i|.
\]
On the other hand, since $\max_i b_i \ge b_{k_0}$, it follows 
\[
a_{k_0}-\max_i b_i\le a_{k_0}-b_{k_0},
\]
which proofs the lemma in Case 1.

\emph{Case 2}: Suppose both $\max_i a_i$ and $\max_i b_i$ are non-positive. We omit the proof of this case as it is similar to the proof of Case 1.

\emph{Case 3}: Suppose $(\max_i a_i)(\max_i b_i)<0$. Without loss of generality, we assume $\max_i a_i>0> \max_i b_i$. Let $1\le k_0\le n$ such that $a_{k_0}=\max_i a_i$. Apparently
\[
|a_{k_0}|+|b_{k_0}|= |a_{k_0}-b_{k_0}|\le \max_{1\le i\le n}|a_i-b_i|.
\]
On the other hand, since $0>\max_i b_i \ge b_{k_0}$, it follows 
\[
\bigg| |a_{k_0}|-\big|\max_i b_i\big| \bigg| \le |a_{k_0}|+\big|\max_i b_i\big| \le  |a_{k_0}|+|b_{k_0}|,
\]
which proofs the lemma in Case 3.

\subsection{Extreme Value Estimation: Proof of Proposition \ref{order.risk.lem}}

Observe that the risk $\E |y_{(p)}-\mu_{(p)}|^\alpha$ is invariant w.r.t. an arbitrary constant mean shift for all $\{\mu_1,...,\mu_p\}$. As a result, without loss of generality, we fix $\mu_{(p)}=0\ge \mu_j$ for all $j=1,...,p,$ and it suffices to evaluate $\E y_{(p)}^\alpha$. Let $y_j=\mu_j+W_j$ where $W_j$ are centred sub-Gaussian random variable with cdf $P_j(w)$ and pdf $p_j(w)$.
By definition, we know that the cdf of the largest order statistic $y_{(p)}$ is
\[
F_{\max}(x)=\prod_{j=1}^pP_j(x-\mu_j),
\]
with the pdf
\[
f_{\max}(x)=\sum_{k=1}^p\prod_{j\ne k}P_j(x-\mu_j)p_k(x-\mu_k).
\]
Thus
\[
\E y_{(p)}^\alpha=\int x^\alpha\sum_{k=1}^p\prod_{j\ne k}P_j(x-\mu_j)p_k(x-\mu_k)dx.
\]
By taking derivative with respect to $\mu_i$, we have for $\mu_i\le 0$,
\begin{align*}
\frac{\partial \E y_{(p)}^\alpha}{\partial \mu_i}&=-\int x^\alpha\sum_{\substack{1\le k\le p,\\ k\ne i}}\prod_{j\ne \{k ,i\}}P_j(x-\mu_j)p_i({x-\mu_i})p_k({x-\mu_k})dx\\
&\quad -\int{x^\alpha}\prod_{j\ne i}P_j({x-\mu_j})p'_i({x-\mu_i})dx\\
&\le 0.
\end{align*}
As a result, the expectation as a function of $(\mu_1,...,\mu_p)$ reaches its maximum only when $\mu_1=\mu_2=...=\mu_p=0$. In this case, defining $Z=y_{(p)}$, it holds that, for $x>0$
\[
P(Z>x)\le pP(|W_1|>x)\le Cp\exp(- x^2/\sigma^2).
\]
Then
\begin{align*}
\E Z^\alpha &= \int_0^\infty \alpha x^{\alpha-1}P(Z>x)dx\\
&\le \int_0^{\sigma\sqrt{c\log p}} \alpha x^{\alpha-1}dx+C\alpha p\int_{\sigma\sqrt{c\log p}}^{\infty}x^{\alpha-1}\exp(- x^2/\sigma^2)dx\\
& = (c\sigma^2\log p)^{\alpha/2}+o(1)
\end{align*}
for some constant $c>0$.
This completes the proof.

\subsection{Calculation of Chi-Square Divergence}

\subsubsection{Proof of Lemma \ref{chisq.lem}}
Let $p_0$ and $p_1$ be the densities corresponding to $P_0$ and $P_1$, respectively. It follows that
\[
p_0=\prod_{i=2}^{n}\prod_{j=1}^p\phi_{0,\sigma}(x_{ij})\cdot \prod_{1\le j\le p}\phi_{\theta_j^{(0)},\sigma}(x_{1j}),\qquad \theta^{(0)}=\bfeta_0
\]
where $\phi_{\mu,\sigma}(x)$ is the density function for the normal distribution $N(\mu,\sigma^2)$ and
\[
p_1=\frac{2}{p}\sum_{k=1}^{p/2}g_k,
\]
with
\[
g_k =\prod_{i=2}^{n}\prod_{j=1}^p\phi_{0,\sigma}(x_{ij})\cdot \prod_{1\le j\le p}\phi_{\theta_j^{(k)},\sigma}(x_{1j}),\qquad \theta^{(k)}=\pi_k(\bfeta_1),
\]
where $\pi_k(\bfeta_1)\in \R^p$ exchanges the last and the $k$-th coordinate of $\bfeta_1$.
Consequently, we have
\[
\int \frac{p_1^2}{p_0} = \frac{4}{p^2} \sum_{1\le i,j \le p/2}\int \frac{g_ig_{j}}{p_0}.
\]
For any $i$ and $j$,
\begin{align*}
\int \frac{g_ig_j}{p_0} &= \frac{1}{\sigma^{p}(2\pi)^{p/2}} \int \exp\bigg\{ -\frac{\sum_{k=1}^p(x_{1k}-\theta_k^{(i)})^2+\sum_{k=1}^p(x_{1k}-\theta_k^{(j)})^2-\sum_{k=1}^p(x_{1k}-\theta_k^{(0)})^2}{2\sigma^2}\bigg\}\\
&\quad \times \frac{1}{\sigma^{p(n-1)}(2\pi)^{\frac{p(n-1)}{2}}} \int \exp\bigg\{ -\frac{\sum_{i=2}^n\sum_{j=1}^px_{ij}^2}{2\sigma^2}\bigg\}\\
&=\frac{1}{\sigma^{p}(2\pi)^{p/2}} \int \exp\bigg\{ -\frac{\sum_{k=1}^p(x_{1k}-\theta_k^{(i)}-\theta_k^{(j)}+\theta_k^{(0)})^2-2\sum_{k=1}^p(\theta_{k}^{(i)}-\theta_k^{(0)})(\theta_k^{(j)}-\theta_k^{(0)})}{2\sigma^2}\bigg\}\\
&=\exp(\log p\cdot I\{i=j \}).
\end{align*}
It then follows that
\[
\int\frac{p_1^2}{p_0} = \frac{4}{p^2} \sum_{1\le i,j\le p/2} \exp(\log p I\{i=j \} ) = \E \exp(\log p\cdot J),
\]
where $J$ is a Bernoulli random variable with $P(J=1)=2/p.$ As a result, we have
\[
\E \exp(J\log p )= (1-2/p)+2e^{\log p}/p \le 3.
\]

\subsubsection{Proof of Lemma \ref{kl.mix.prop}}

Again, the direct calculation of $D(\bar{P}_{u,t,\beta_R}\| \bar{P}_{u,t,\beta_R})$ is difficult, so we detour by introducing an approximate density of $\bar{P}_{u,t,\beta_R}$ as
\begin{align*}
\tilde{P}_{u,t,\beta_R}(Y)=\frac{1}{(2\pi)^{np/2}\sigma^{np}}\int&\exp(-\|Y-2tuw_{+}^\top\|_F^2/(2\sigma^2))\\
&\times \bigg( \frac{p-1}{2\pi}\bigg)^{(p-1)/2}\exp(-(p-1)\|w\|_2^2/2)dw.
\end{align*}
Let $Y=[Y_1\quad Y_2]$ where $Y_1\in \R^{n\times (p-1)}$ and $Y_2\in \R^{n\times 1}$.The above expression can be simplified in the sense that, for $Y\sim \tilde{P}_{u,t,\beta_R }$, if $Y_i$ is the $i$-th column of $Y\in \R^{n\times p}$, we have
\beq \label{Y_i.dist}
Y_i| u \sim_{i.i.d.} N\bigg( 0, \sigma^2\bigg(I_n-\frac{4t^2}{4t^2+\sigma^2(p-1)}u u^\top\bigg)^{-1} \bigg) = N\bigg(0,\sigma^2I_n+\frac{4t^2}{(p-1)}uu^\top\bigg), 
\eeq
for $i=1,...,p-1,$ and
\[
Y_p|u \sim N(2t\beta_R u, \sigma^2I_n)\quad \text{independent of $Y_i, i\ne p$.}
\]
It is well-known that the KL-divergence between two $p$-dimensional multivariate Gaussian distribution is
\[
D( N(\mu_0,\Sig_0)\| N(\mu_1,\Sig_1)) =\frac{1}{2}\bigg(  \text{tr}(\Sig_0^{-1}\Sig_1)+(\mu_1-\mu_0)^\top \Sig_1^{-1}(\mu_1-\mu_0)-p+\log\bigg( \frac{\det \Sig_1}{\det \Sig_0}\bigg) \bigg).
\]
As a result, we can calculate that for any $\tilde{P}_{u,t,\beta_R }$ and $\tilde{P}_{u',t,\beta_R }$,
\begin{align} \label{KL.tilde}
D(\tilde{P}_{u,t,v_p}\| \tilde{P}_{u',t,v_p}) &=\frac{p-1}{2}\bigg\{ \text{tr}\bigg( \bigg( I_n -\frac{4t^2}{4t^2+\sigma^2(p-1)}u u^\top \bigg)\bigg( I_n +\frac{4t^2}{\sigma^2(p-1)}u' {u'}^\top \bigg)\bigg)-n \bigg\}\nonumber \\
&\quad +\frac{ Ct^2\beta_R^2\|u-u'\|_2^2}{\sigma^2}\nonumber \\
&\le \frac{Ct^4}{(4t^2+\sigma^2(p-1))}(1-(u^\top u')^2),
\end{align}
where the last inequality follows from $t^2/(t^2+\sigma^2p)\ge \beta_R^2$.
Hence, the proof of this proposition is complete if we can show that there exist some constant $C>0$ such that
\beq \label{affinity.2}
D(\bar{P}_{u,t,\beta_R}\| \bar{P}_{u,t,\beta_R})\le D(\tilde{P}_{u,t,\beta_R}\| \tilde{P}_{u,t,\beta_R})+C.
\eeq
The rest of the proof is devoted to the proof of (\ref{affinity.2}).
\paragraph{Proof of (\ref{affinity.2}).} For any given $u$,
\begin{align} \label{pu/qu}
&\frac{\bar{P}_{u,t,\beta_R}}{\tilde{P}_{u,t,\beta_R}}\\
&=\frac{1}{(2\pi)^{\frac{p-1}{2}}(\frac{\sigma^2}{4t^2+\sigma^2(p-1)})^{\frac{p-1}{2}}}\exp\bigg(\frac{1}{2\sigma^2}\|Y_2-6t\beta_R u \|_F^2+\frac{1}{2\sigma^2} \sum_{i=1}^{p-1}Y_{1,i}^\top(I_n-\frac{4t^2}{4t^2+\sigma^2(p-1)}uu^\top) Y_{1,i}\bigg) \nonumber \\
&\quad\times C_{t,\beta_R } \int_{\mathcal{G}}\exp(-\|Y-2tuw_{+}^\top\|_F^2/(2\sigma^2)-(p-1)\|w\|_F^2/2)dw \nonumber \\
&=\frac{1}{(2\pi)^{(p-1)/2}(\frac{\sigma^2}{4t^2+\sigma^2(p-1)})^{(p-1)/2}}\exp\bigg(\frac{1}{2\sigma^2} \sum_{i=1}^{p-1}Y_{1,i}^\top(I_n-\frac{4t^2}{4t^2+\sigma^2(p-1)}uu^\top) Y_{1,i}\bigg)\nonumber \\
&\quad\times C_{t,\beta_R } \int_{\mathcal{G}}\exp\bigg(-\frac{1}{2\sigma^2}\bigg[\text{tr}(Y_1^\top(I_n-\frac{4t^2}{4t^2+\sigma^2(p-1)}uu^\top) Y_1\nonumber \\
&\quad+(4t^2+\sigma^2(p-1))(w-\frac{2t}{4t^2+\sigma^2(p-1)}Y_1^\top u)(w-\frac{2t}{4t^2+\sigma^2(p-1)}Y_1^\top u)^\top)\bigg]\bigg)dw\nonumber \\
&=C_{t,\beta_R }P\bigg(w\in\mathcal{G}\bigg|  w \sim N\bigg( \frac{2t}{4t^2+\sigma^2(p-1)}Y_1^\top u, \frac{\sigma^2}{4t^2+\sigma^2(p-1)}{\bf I}_{n} \bigg)\bigg)\nonumber \\
&\le C_{t,\beta_R }.
\end{align}
where we recall that $Y_1$ is defined by $Y=[Y_1\quad Y_2]$ with $Y_1\in \R^{n\times (p-1)}$ and $Y_2\in \R^{n\times 1}$. Recall that
\[
C^{-1}_{t,\beta_R }= P\big( w\in\mathcal{G}| w_j\sim N(0,1/(p-1))  \big).
\]
By concentration inequality for sub-exponential random variables, we have
\[
P\bigg( \bigg|\frac{p-1}{p} \|w\|_2^2-1\bigg|\le \sqrt{\frac{\log p}{p}} \bigg)\ge 1-O(p^{-c}),
\]
and therefore for sufficiently large $p$,
\[
P(1/4\le \|w\|_2^2\le 4) \ge P\bigg(  1- \sqrt{\frac{\log p}{p}} \le  \|w\|_2^2\le 1+\sqrt{\frac{\log p}{p}} \bigg)\ge 1-O(p^{-c}).
\]
Moreover, by standard Gaussian tail bounds we have $|\bar{w}|\lesssim \sqrt{\log p}/p$ and $w_{(p)} \lesssim p^{-1/2}\sqrt{\log p}$ with probability at least $1-O(p^{-c})$.
Then it follows that 
\beq \label{c-1}
C^{-1}_{t,\beta_R }\ge 1-p^{-c}
\eeq
for some small constant $c>0$.
By (\ref{c-1}), we know that
\[
\frac{\bar{P}_{u,t,\beta_R}}{\tilde{P}_{u,t,\beta_R}}\le 1+p^{-c}
\]
uniformly for some constant $c>0.$ Thus, for some constant $\delta>0$, we have
\begin{align*}
&D(\bar{P}_{u,t,\beta_R}\| \bar{P}_{u',t,\beta_R})\\
&=\int \bar{P}_{u,t,\beta_R}\bigg[\log \bigg( \frac{\bar{P}_{u,t,\beta_R}}{\tilde{P}_{u,t,\beta_R}} \bigg)+\log \bigg( \frac{\tilde{P}_{u,t,\beta_R}}{\tilde{P}_{u',t,\beta_R}} \bigg)+\log \bigg( \frac{\tilde{P}_{u',t,\beta_R}}{\bar{P}_{u',t,\beta_R}} \bigg)\bigg]dY\\
&\le \log(1+\delta)+D(\tilde{P}_{u,t,\beta_R}\| \tilde{P}_{u',t,\beta_R}) +\int (\bar{P}_{u,t,\beta_R}-\tilde{P}_{u,t,\beta_R})\log \bigg( \frac{\tilde{P}_{u,t,\beta_R}}{\tilde{P}_{u',t,\beta_R}} \bigg)dY+\int\bar{P}_{u,t,\beta_R}\log \bigg( \frac{\tilde{P}_{u',t,\beta_R}}{\bar{P}_{u',t,\beta_R}} \bigg)dY \\
&\le \log(1+\delta)+D(\tilde{P}_{u,t,\beta_R}\| \tilde{P}_{u',t,\beta_R}) +\int\tilde{P}_{u,t,\beta_R} \bigg(\frac{\bar{P}_{u,t,\beta_R}}{\tilde{P}_{u,t,\beta_R}}-1\bigg)\log \bigg( \frac{\tilde{P}_{u,t,\beta_R}}{\tilde{P}_{u',t,\beta_R}} \bigg)dY\\
&\quad+(1+\delta)\int\tilde{P}_{u,t,\beta_R}\bigg|\log \bigg( \frac{\tilde{P}_{u',t,\beta_R}}{\bar{P}_{u',t,\beta_R}} \bigg)\bigg|dY\\
&\le \log(1+\delta)+D(\tilde{P}_{u,t,\beta_R}\| \tilde{P}_{u',t,\beta_R}) +p^{-c}\int\tilde{P}_{u,t,\beta_R} \bigg|\log \bigg( \frac{\tilde{P}_{u,t,\beta_R}}{\tilde{P}_{u',t,\beta_R}} \bigg)\bigg|dY\\
&\quad+(1+\delta)\int\tilde{P}_{u,t,\beta_R}\bigg|\log \bigg( \frac{\tilde{P}_{u',t,\beta_R}}{\bar{P}_{u',t,\beta_R}} \bigg)\bigg|dY.
\end{align*}
Now since
\begin{align*}
&\int\tilde{P}_{u,t,\beta_R} \bigg|\log \bigg( \frac{\tilde{P}_{u,t,\beta_R}}{\tilde{P}_{u',t,\beta_R}} \bigg)\bigg|dY\\
&=\frac{1}{2\sigma^2}\int\tilde{P}_{u,t,\beta_R} \bigg|\frac{4t^2}{4t^2+\sigma^2(p-1)}\sum_{i=1}^{p-1}Y_{.i}^\top(uu^\top-u'{u'}^\top)Y_{.i}\bigg|dY\\
&\le \frac{1}{2\sigma^2}\E\bigg[ \frac{4t^2}{4t^2+\sigma^2(p-1)}\sum_{i=1}^{p-1}Y_{.i}^\top(uu^\top+u'{u'}^\top)Y_{.i} \bigg|Y_{.i}\sim N\bigg(0,\sigma^2I_n+\frac{4t^2\sigma^2}{\sigma^2(p-1)}uu^\top \bigg)\bigg]\\
&=\frac{4t^2(p-1)}{2\sigma^2(4t^2+\sigma^2(p-1))}\text{tr}\bigg( (uu^\top+u'{u'}^\top)\big(\sigma^2I_n+\frac{4t^2}{(p-1)}uu^\top\big)  \bigg)\\
&\le \frac{4t^2(p-1)}{4t^2+\sigma^2(p-1)}\text{tr}\bigg(u^\top\big(I_n+\frac{4t^2}{\sigma^2(p-1)} \big)u  \bigg)\\
&=\frac{4t^2}{\sigma^2}\le  p,
\end{align*}
we know that the third term in the above expression can be bounded by
\[
p^{-c}\int\tilde{P}_{u,t,\beta_R} \bigg|\log \bigg( \frac{\tilde{P}_{u,t,\beta_R}}{\tilde{P}_{u',t,\beta_R}} \bigg)\bigg|dY\le p\cdot p^{-c}\le C
\]
for some constant $C>0$ and sufficiently large $c>0$.
Finally, by (\ref{pu/qu}), we have
\begin{align*}
\int \tilde{P}_{u,t,\beta_R}\bigg|\log \bigg( \frac{\tilde{P}_{u',t,\beta_R}}{\bar{P}_{u',t,\beta_R}} \bigg)\bigg|dY&\le \int\tilde{P}_{u,t,\beta_R}\bigg|\log \frac{1}{C_{t,\beta_R}}\bigg|dY+\int \tilde{P}_{u,t,\beta_R}\bigg|\log \frac{1}{P(w\in\mathcal{G}|E)}\bigg|dY,
\end{align*}
where we denoted 
\[
E=\bigg\{ w \sim N\bigg( \frac{2t}{4t^2+\sigma^2(p-1)}Y_1^\top u', \frac{\sigma^2}{4t^2+\sigma^2(p-1)}{\bf I}_{n} \bigg) \bigg\}.
\]
Now on the one hand,
\[
\int\tilde{P}_{u,t,\beta_R}\bigg|\log \frac{1}{C_{t,\beta_R}}\bigg|dY\le \big(\log(1+\delta) \lor |\log(1-\delta)^{-1}|\big).
\]
On the other hand, if we denote $\mathcal{G}_1=\{ 1/2\le \|w\|_2 \}$, $\mathcal{G}_2=\{ \|w\|_2\le 2 \}$, $\mathcal{G}_3=\{  w_{(p)}\le cp^{-1/2}\sqrt{\log p} \}$ and $\mathcal{G}_4=\{ \bar{w}\le cp^{-1/2}\sqrt{\log p}\}$, and denote $P(A_i)=P(w\in \mathcal{G}_i|E)$ for $i=1,...,4$, then by definition, $\mathcal{G}=\bigcap_{i=1}^4\mathcal{G}_i$, and direct calculation yields
\begin{align}
P(w\in\mathcal{G}| E)&=P(A_1,A_2|A_3,A_4)P(A_4|A_3)P(A_3)\nonumber \\
&\ge P(A_3)P(A_1,A_2|A_3)\qquad [A_3\subseteq A_4] \nonumber \\
&\ge P(A_3)P(A_1,A_2)\qquad [A_1\cap A_2 \subset A_3]. \label{GG}
\end{align}
Hence,
\beq \label{last.b}
\int \tilde{P}_{u,t,\beta_R}\bigg|\log \frac{1}{P(w\in\mathcal{G}|E)}\bigg|dY\le -\int \tilde{P}_{u,t,\beta_R}\log \frac{1}{P(A_3)}dY- \int \tilde{P}_{u,t,\beta_R}\log \frac{1}{P(A_1,A_2)}dY
\eeq
For the first term in (\ref{last.b}), we consider the lower bound for the probability
\[
P(A_3)=P(w_{(p)}\le cp^{-1/2}\sqrt{\log p}|E )\ge P(w_{(p)}\sqrt{4t^2+\sigma^2 p}\le \sigma c\sqrt{\log p}|E ).
\]
Denote $\Delta(Y_1)=c\sqrt{\log p}-\frac{2t\|Y_1^\top u'\|_\infty}{\sigma\sqrt{4t^2+\sigma^2p}}$. We have
\begin{align*}
&P\bigg(\sqrt{4t^2+\sigma^2p}w_{(p)}/\sigma \le c\sqrt{\log p}\bigg|E, \Delta(Y_1)\ge \sqrt{\log p} \bigg)\\
&\ge 1-p\max_iP\bigg(\sqrt{4t^2+\sigma^2p}w_{i}/\sigma > c\sqrt{\log p}\bigg|E, \Delta(Y_1)\ge \sqrt{\log p} \bigg)\\
&\ge 1-p\Phi\bigg(- c\sqrt{\log p}+\frac{2t\|Y_1^\top u'\|_\infty}{\sqrt{4t^2+\sigma^2(p-1)}} \bigg)\\
&\ge 1/2
\end{align*}
for sufficiently large $p>0$.
Moreover,
\begin{align*}
&P\bigg(\max_j\sqrt{4t^2+\sigma^2p}w_{j}/\sigma \le c\sqrt{\log p}\bigg|E,\Delta(Y_1)\le \sqrt{\log p} \bigg) \\
&=\Phi^p\bigg( c\sqrt{\log p}-\frac{t(Y_1^\top u')_j}{\sigma\sqrt{4t^2+\sigma^2p}} \bigg)\\
&\ge \Phi^p\bigg( c\sqrt{\log p}-\frac{t\|Y_1^\top u'\|_\infty}{\sigma\sqrt{4t^2+\sigma^2p}} \bigg)\\
&\gtrsim e^{-p\Delta^2}.
\end{align*}
As a result
\beq \label{G3}
P(A_3)\ge  0.5\cdot I\{  \Delta(Y_1)> \sqrt{\log p}\}+C\exp(-p\Delta^2)I\{  \Delta(Y_1)\le \sqrt{\log p} \},
\eeq
and we have
\begin{align} \label{pa_3}
-\int \tilde{P}_{u,t,\beta_R}\log \frac{1}{P(A_3)}dY\le -\log(1/2)\int_{\Delta(Y_1)> \sqrt{\log p}} \tilde{P}_{u,t,\beta_R}dY+Cp\int_{\Delta(Y_1)\le \sqrt{\log p}}\tilde{P}_{u,t,\beta_R} \Delta^2(Y_1) dY.
\end{align}
where the first term is obviously bounded. Note that $\Delta(Y_1)\le \sqrt{\log p}$ implies
\[
\|Y_1^\top u'\|_\infty \ge (c+1)\sqrt{\log p}\cdot \frac{\sigma\sqrt{4t^2+\sigma^2p}}{2t}.
\]
The second term can bee bounded by
\begin{align*}
p\int_{\Delta(Y_1)\le \sqrt{\log p}}\tilde{P}_{u,t,\beta_R} \Delta^2(Y_1) dY&\le p\int_{\frac{2t\|Y_1^\top u'\|_\infty}{\sqrt{4t^2+\sigma^2 p}} \ge \sigma(c+1)\sqrt{\log p}}\tilde{P}_{u,t,\beta_R} \Delta^2(Y_1) dY \\
&\le p\int_{\frac{2t\|Y_1^\top u'\|_\infty}{\sqrt{4t^2+\sigma^2 p}} \ge \sigma(c+1)\sqrt{\log p}}\tilde{P}_{u,t,\beta_R} \big(\log p\lor \frac{4t^2\|Y_1^\top u'\|_\infty^2}{\sigma^2(4t^2+\sigma^2p)}\big)dY\\
&\le p\int_{\frac{2t\|Y_1^\top u'\|_\infty}{\sqrt{4t^2+\sigma^2 p}} \ge \sigma(c+1)\sqrt{\log p}}\tilde{P}_{u,t,\beta_R}  \frac{4t^2\|Y_1^\top u'\|_\infty^2}{\sigma^2(4t^2+\sigma^2p)}dY\\
&=p\E \|\xi\|_\infty^2I\{\|\xi\|_\infty\ge (c+1)\sqrt{\log p} \}
\end{align*}
where in the last inequality, 
\[
\xi\sim N\bigg(0,\frac{4t^2(\sigma^2(p-1)+4t^2(u^\top u')^2)}{\sigma^2 (p-1)(\sigma^2p+4t^2)}I_{p-1}\bigg),
\]
which follows from (\ref{Y_i.dist}) and the fact that
\[
Y_1^\top u\sim N\bigg(0, \frac{\sigma^2 (p-1)+4t^2(u^\top u')^2}{p-1}I_{p-1}\bigg).
\]
Now since
\[
P(\|\xi\|_\infty>x)\le pP(|\xi_i|>x)=2p\Phi(-x/\tau)<2p \tau x^{-1}\phi(x/\tau),
\]
\[
\tau^2=\frac{4t^2(\sigma^2(p-1)+4t^2(u^\top u')^2)}{\sigma^2 (p-1)(\sigma^2p+4t^2)},
\]
we have
\begin{align*}
\E \|\xi\|_\infty^2I\{\|\xi\|_\infty\ge (c+1)\sqrt{\log p} \} &= \int_{(c+1)\sqrt{\log p}}^\infty 2 xP(\|\xi\|_\infty>x)dx\\
&\le 4p\tau\int_{(c+1)\sqrt{\log p}}^{\infty}\phi(x/\tau)dx\\
& \lesssim \tau^2p^{-\frac{(c+1)^2}{\tau^2}+1}.
\end{align*}
Note that by definition $\tau\le 1$. As long as $c\ge 1$, we have
\[
\frac{(c+1)^2}{\tau^2}\ge 2,
\]
and therefore
\[
p\int_{\Delta(Y_1)\le \sqrt{\log p}}\tilde{P}_{u,t,\beta_R} \Delta^2(Y_1) dY=O(1).
\]
Then (\ref{pa_3}) implies
\[
-\int \tilde{P}_{u,t,\beta_R}\log \frac{1}{P(A_3)}dY =O(1).
\]
Lastly, we considerthe second term in (\ref{last.b}). For the event $A_1$, since $Y_1^\top u\in \R^{(p-1)\times 1}$, we can find a orthogonal matrix $Q\in \R^{(p-1)\times (p-2)}$ such that $Q^\top Y_1^\top u=0$ and hence $Q^\top w\sim N(0,\frac{\sigma^2}{4t^2+\sigma^2(p-1)})$. Then by standard result in random matrix (e.g. Corollary 5.35 in \cite{vershynin2010introduction} and treating $Q^\top w$ as a matrix), we have
\[
\|w\|_2\ge \| Q^\top w\|_{2} \ge \frac{\sigma}{\sqrt{4t^2+\sigma^2(p-1)}}(\sqrt{p-2}-1-x)
\] 
with probability at least $1-2\exp(-x^2/2)$. Since $t^2<\sigma^2 p/4$, for $p$ sufficiently large, we can find $c$ such that by setting $x=c\sqrt{p}$, 
\beq \label{G1}
P( A_1) \ge 1-e^{-cp}.
\eeq 
For the event $A_2$, denote
\[
\Gamma_s(Y_1)=\frac{s(4t^2+\sigma^2 (p-1))}{\sigma^2}-\frac{4t^2\|Y_1^\top u\|_2^2}{\sigma^2(4t^2+\sigma^2 (p-1))}.
\]
If $\Gamma_4(Y_1)\ge p$, since
by concentration inequality for sub-exponential random variables, we have for any $x>0$,
\[
P\bigg( \bigg|\|w\|_2^2 -\frac{4t^2\|Y_1^\top u\|_2^2}{(4t^2+\sigma^2(p-1))^2} \bigg| \le \frac{\sigma^2x}{4t^2+\sigma^2(p-1)}\bigg)\ge 1-e^{-c (x^2/n\land x)},
\]
or
\[
P\bigg( \|w\|_2^2 \ge \frac{4t^2\|Y_1^\top u\|_2^2}{(4t^2+\sigma^2(p-1))^2} + \frac{\sigma^2x}{4t^2+\sigma^2(p-1)}\bigg)\le e^{-c (x^2/n\land x)}.
\]
Then by setting $x=\Gamma_4(Y_1)$, we have
\[
P( \|w\|_2^2 \ge 4| \Gamma_4(Y_1)\ge p)\le e^{-c (p^2/n\land p)}.
\]
If $\Gamma_4(Y_1)< p$, since by the tail upper and lower bound for sum of independent sub-exponential random variables (\cite{zhang2018non}), we have for any $x>0$
\beq \label{tail}
e^{-c_1x}\ge P\bigg( \|w\|_2^2 \le \frac{4t^2\|Y_1^\top u\|_2^2}{(4t^2+\sigma^2(p-1))^2} - \frac{\sigma^2x}{4t^2+\sigma^2(p-1)} \bigg)\ge e^{-c_2 (x^2/n\land x)}\ge e^{-c_2x}.
\eeq
Set $x=-\Gamma_4(Y_1)$. If in addition, $\Gamma_4(Y_1)<0$, then by the above tail probability lower bound,
\[
P\bigg( \|w\|_2^2 \le 4\bigg| \Gamma_4(Y_1)<0 \bigg)\ge e^{c_2 \Gamma_4}.
\]
If $0\le \Gamma_4(Y_1)<p$, then
\[
P\bigg( \|w\|_2^2 \le 4 \bigg| 0\le \Gamma_4(Y_1)<p\bigg)\ge 1/2.
\]
By (\ref{tail}), we have
\begin{align*}
P\bigg( 1/4\le \|w\|_2^2 \le 4 \bigg| \Gamma_4(Y_1)<0  \bigg)&=P\bigg( \|w\|_2^2 \le 4 \bigg|\Gamma_4(Y_1)<0 \bigg)-P\bigg(  \|w\|_2^2\le 1/4 \bigg|\Gamma_4(Y_1)<0 \bigg)\\
&\ge  \exp\{c_2 \Gamma_4(Y_1)\}- \exp\{c_1 \Gamma_{1/4}(Y_1)\}.
\end{align*}
Note that the choice of the constant upper and lower bound $(a,b)$ (chosen as $(1/4,4)$ here) are not essential throughout the proof. Basically, they should depend on the constant $c_1,c_2$ in the last expression so that
\beq \label{pos}
\exp\{c_2 \Gamma_b(Y_1)\}- \exp\{c_1 \Gamma_{1}(Y_1)\}\ge 0,\qquad bc_2>ac_1.
\eeq
Without loss of generality, we can assume that $(1/4,4)$ is an admissible pair for (\ref{pos}) to be true. As a result,
\[
\exp\{c_2 \Gamma_4(Y_1)\}- \exp\{c_1 \Gamma_{1/4}(Y_1)\}\gtrsim \exp\bigg\{- \frac{4c_2t^2\|Y_1^\top u\|_2^2}{\sigma^2(t^2+\sigma^2p)} \bigg\} \cdot pe^{cp}
\]
Thus,
\begin{align*}
&-\int \tilde{P}_{u,t,\beta_R}\log P(A_1,A_2)dY\\
&=-\int \tilde{P}_{u,t,\beta_R}\log [P(A_1)-P(A_2^c)]dY\\
&\le -\int_{\Gamma_4\ge p} \tilde{P}_{u,t,\beta_R}\log [1-e^{-cp}-P(A_2^c)]dY-\int_{\Gamma_4< p} \tilde{P}_{u,t,\beta_R}\log [P(A_1,A_2)]dY\\
&\le -\log(1-2e^{-cp})-\int_{\Gamma_4< 0} \tilde{P}_{u,t,\beta_R}\log [P(1/2\le \|w\|_2\le 2|E)]dY\\
&\quad-\int_{0\le \Gamma_4<p} \tilde{P}_{u,t,\beta_R}\log [P(1/2\le \|w\|_2\le 2|E)]dY,
\end{align*}
where
\begin{align*}
-\int_{0\le \Gamma_4<p} \tilde{P}_{u,t,\beta_R}\log [P(1/2\le \|w\|_2\le 2|E)]dY&=-\int_{0\le \Gamma_4<p} \tilde{P}_{u,t,\beta_R}\log [P(A_2)-P(A_1^c)]dY\\
&\le -\log(1/2-e^{-cp})
\end{align*}
and
\begin{align*}
&-\int_{\Gamma_4< 0} \tilde{P}_{u,t,\beta_R}\log [P(1/2\le \|w\|_2\le 2|E)]dY\\
&\le-\int_{\Gamma_4< 0} \tilde{P}_{u,t,\beta_R}\log \bigg(  \exp\bigg\{- \frac{4c_2t^2\|Y_1^\top u\|_2^2}{\sigma^2(t^2+\sigma^2p)} \bigg\} \cdot pe^{cp}\bigg)dY\\
&\le c_2\int_{\Gamma_4< 0} \tilde{P}_{u,t,\beta_R}\frac{4t^2\|Y_1^\top u\|_2^2}{\sigma^2(t^2+\sigma^2p)} dY+cP(\Gamma_4< 0 )\\
&\lesssim \E \|\xi\|_2^2I\{ \Gamma_4< 0 \}+c.
\end{align*}
Finally, we calculate
\[
\E \|\xi\|_2^2I\{\Gamma_4< 0  \}\le \int_{8 p}^\infty P(\|\xi\|_2^2\ge x)dx=O(1).
\]
So the second term in (\ref{last.b}) is also bounded by some constant. In sum, we have proven the inequality (\ref{affinity.2}) and therefore completed the proof.

\section{Supplementary Tables and Figures}

\begin{figure}[h!]
	\centering
	\includegraphics[angle=0,width=4.5cm]{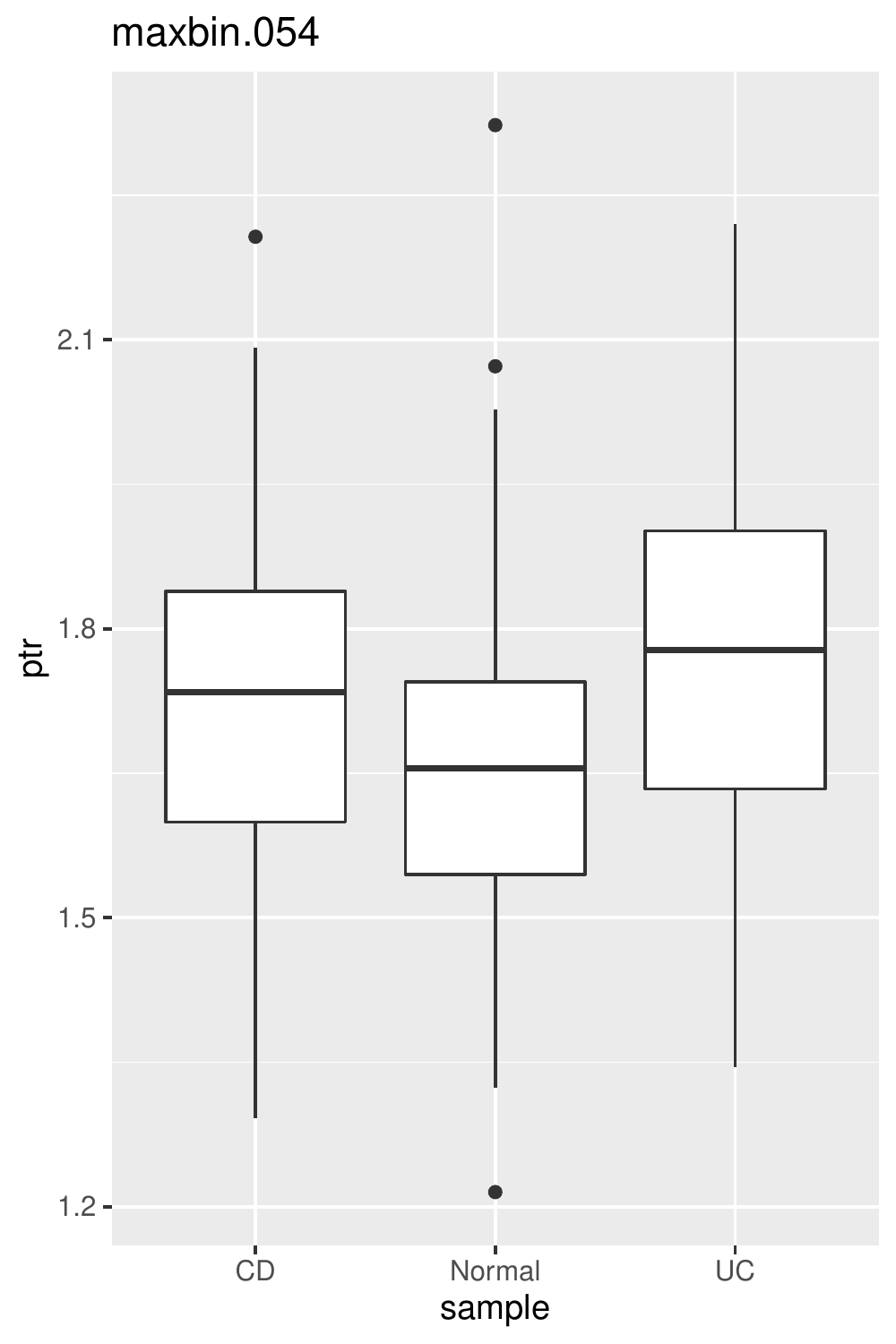}
	\includegraphics[angle=0,width=4.5cm]{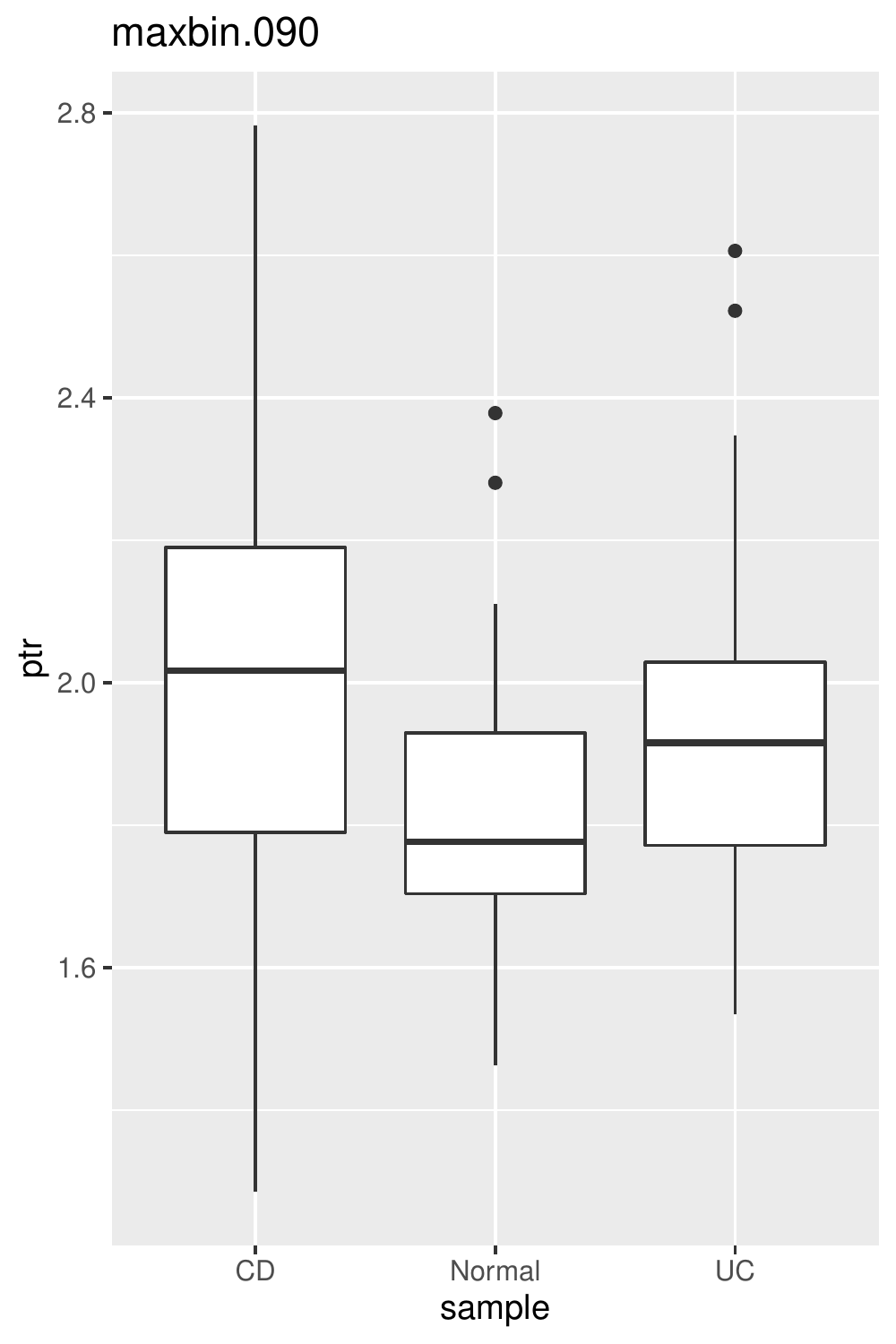}
	\includegraphics[angle=0,width=4.5cm]{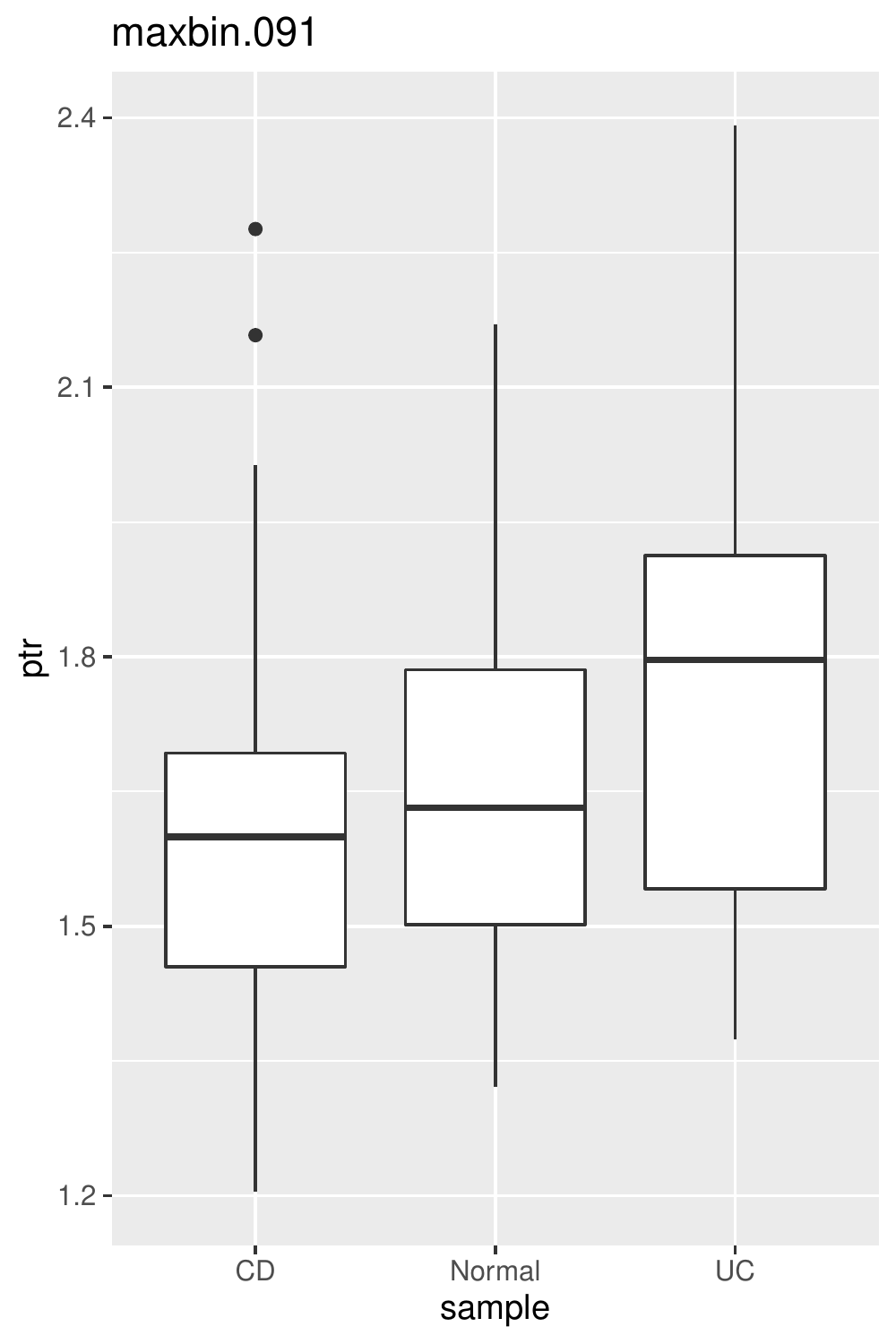}
	\includegraphics[angle=0,width=4.5cm]{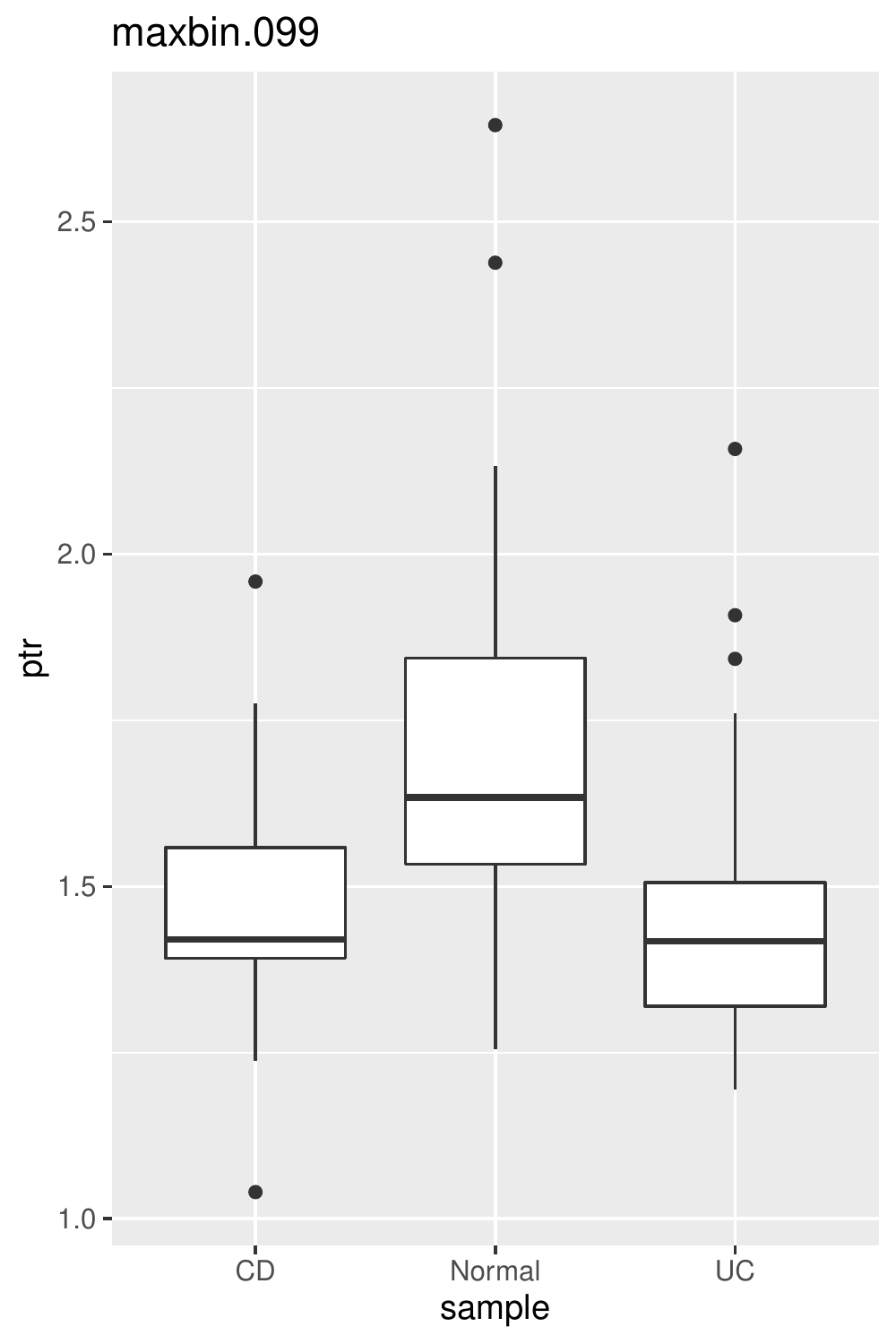}
	\includegraphics[angle=0,width=4.5cm]{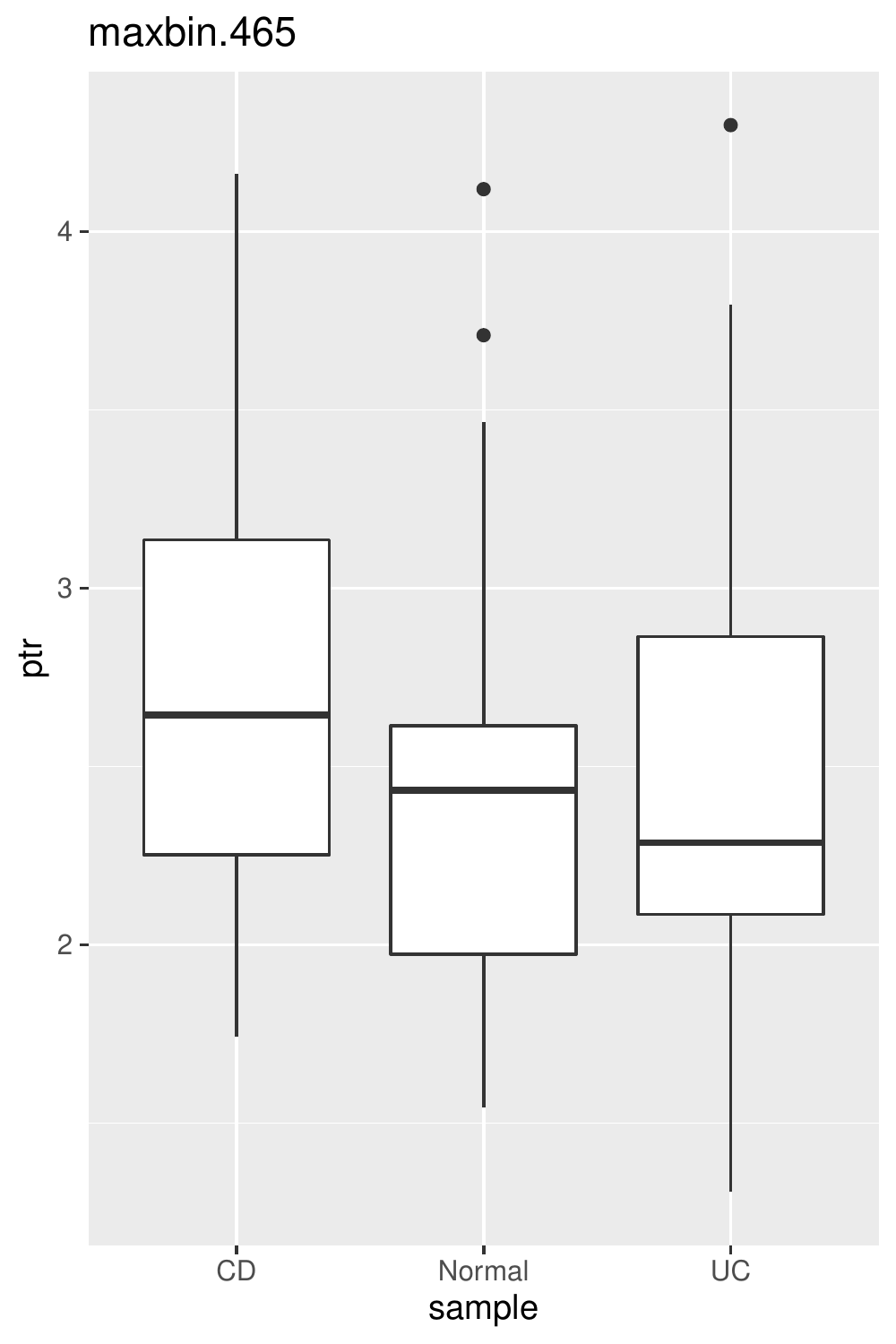}
	\caption{Boxplots of the ePTRs of the five differential bins identified from the F test} 
	\label{simu.3}
\end{figure}

We summarize our supplementary tables and figures as follows:
\begin{enumerate}
	\item In Table \ref{tax.ann}, we provide the complete taxonomic annotations of the five differential bins identified from the F test in Section 5.3 of the main paper.
	\item In Figure \ref{simu.3}, we provide the boxplots of the ePTRs of the above five differential bins.
\end{enumerate}

\begin{table}[h!]
	\centering
	\renewcommand\thetable{S.1}
	\caption{The complete taxonomic annotations with lineage scores indicating the quality of each taxonomic classification.}
	\vspace{0.3cm}
	\begin{tabular*}{ 1.0 \textwidth}{ll}
		\hline 
		Bins  &  Taxonomic Annotations with Lineage Scores  \\  
		\hline  
		bin.054 & Firmicutes (phylum): 0.98;   Clostridia (class): 0.97;  Clostridiales (order): 0.97\\	&Lachnospiraceae (family): 0.93;  Roseburia (genus): 0.92\\
		bin.090 & Firmicutes (phylum): 0.97;	Clostridia (class): 0.95;	Clostridiales (order): 0.95;\\	&Ruminococcaceae (family): 0.84;	Faecalibacterium (genus): 0.82\\
		bin.091  &  Firmicutes (phylum): 0.95; Clostridia (class): 0.93; Clostridiales (order): 0.93\\
		bin.099 & Firmicutes (phylum): 0.96; Clostridia (class): 0.94; Clostridiales (order): 0.94;\\	&Ruminococcaceae (family): 0.87; Subdoligranulum (genus): 0.76;\\
		&	Subdoligranulum sp. APC924/74 (species): 0.74\\
		bin.465 & Firmicutes (phylum): 0.98; Negativicutes (class): 0.97; Veillonellales (order): 0.96;\\	&Veillonellaceae (family): 0.96;  Dialister (genus): 0.94\\
		\hline
	\end{tabular*}
	\label{tax.ann}
\end{table}

\label{lastpage}

\end{document}